\documentclass[12pt,reqno]{amsart}

\usepackage{amssymb,amsmath,graphicx,amsfonts,euscript,listings}
\usepackage{color}
\usepackage{cite}

\allowdisplaybreaks

\newcommand{\Om}{\Omega}

\newcommand{\na}{\nabla}

\def\B{\tilde{b}}
\def\U{\tilde{u}}
\newcommand{\p}{\partial}
\newcommand{\C}{\cdot}

\def\De{\Delta}
\def\na{\nabla}

\newcommand{\ddt}{\frac{d}{dt}}

\newtheorem{thm}{Theorem}[section]
\newtheorem{prop}[thm]{Proposition}
\newtheorem{lem}[thm]{Lemma}

\newtheorem{re}[thm]{Remark}
\newtheorem{no}[thm]{Notation}

\newcommand{\beq}{\begin{equation}}
\newcommand{\eeq}{\end{equation}}
\newcommand{\ben}{\begin{eqnarray}}
\newcommand{\een}{\end{eqnarray}}
\newcommand{\beno}{\begin{eqnarray*}}
\newcommand{\eeno}{\end{eqnarray*}}
\everymath{\displaystyle}

\makeatletter
\@namedef{subjclassname@2020}{%
	\textup{2020} Mathematics Subject Classification}
\makeatother

\numberwithin{equation}{section}

\begin{document}

\title{ Global well-posedness of the 3D Primitive Equations with magnetic field}
\author[Lili Du and  Dan Li]{Lili Du$^{1}$ and  Dan Li$^{2}$}

\address{$^1$ Department of Mathematics, Sichuan University, Chengdu 610064, P.R. China}

\email{dulili@scu.edu.cn}

\address{$^2$ Department of Mathematics, Sichuan University, Chengdu 610064, P.R. China}

\email{dandy0219@hotmail.com}

\vskip .2in
\begin{abstract}
 In this paper, the three-dimensional primitive equations with magnetic field (PEM) are considered on a thin domain. We showed the global existence and uniqueness (regularity) of strong solutions to the  three-dimensional incompressible PEM without any small assumption on the initial data. More precisely, there exists a unique strong solution globally in time for any given $H^2$-smooth initial data.
\end{abstract}
\maketitle

\section{Introduction and the main results}
\subsection{Background and motivation}
In the context of the geophysical flow concerning the large-scale oceanic dynamics, the ratio of the depth to the horizontal width is very small. With the aid of this fact, by scaling the incompressible Navier-Stokes equations with respect to the aspect ratio parameter and taking the small aspect ratio limit, one obtains formally the primitive equations for the large-scale oceanic dynamics. The rigorous mathematical justification of the small aspect ratio limit from the Navier-Stokes equations converges to the primitive equations, which was studied by  Az\'erad-Guill\'en in \cite{AF}. By relying on the result in \cite{AF} to prove the weak convergence, it was  shown in \cite{JKL} that the Navier-Stokes equations strongly converge to the primitive equations. The primitive equations are widely considered as the basic equations of atmospheric dynamics in meteorology. These equations  are the foundation in the weather prediction models, see \cite{GR, RL, AM, JP, GK, WM, QC}. The mathematical analysis of primitive equations was initialed in 1990s by Lions, Temam, and Wang in  \cite{JR, JL, JL1}, where  they established the global existence of weak solutions. The uniqueness of weak solutions for the two-dimensional case was later proved by Bresch et al. in \cite{DBFG}. However, the uniqueness of  the weak solution for the three-dimensional case is still unclear. An important progress for the global well-posedness of  the strong solutions to the three-dimensional primitive equations in a general cylindrical domain has been made by Cao and Titi in \cite{CT}.  This observation and careful study  for the primitive equations allow to establish the well-posedness theory with different boundary conditions and partial viscosity and diffusivity in \cite{JKLT, CLT1, CLT2, KI} and references therein for various generalizations.

\vskip .1in

Compared with the Navier-Stokes equations, the magnetohydrodynamic (MHD) equations, which consist of the Navier-Stokes equations of fluid dynamics and Maxwell's equations of electromagnetism, contain much richer structure than Navier-Stokes equations. The MHD equations reflect the basic physics laws governing the motion of electrically conducting fluids such as
plasmas, liquid metals, and electrolytes, and have played pivotal roles in the study of many phenomena in geophysics, astrophysics, cosmology and engineering (see, e.g., \cite{DB, PAD}). Beside their wide physical applicability, the MHD equations are also of great interest in mathematics. The existence and uniqueness results for weak and strong solutions of the 2D MHD equations are well known in Duvaut and Lions \cite{DL}. For 3D case, it is currently unknown whether the solutions  can develop finite time singularities even if the initial value is sufficiently smooth. Different criteria for regularity in terms of the velocity field, the magnetic field, the pressure or their derivatives have been proposed (see \cite{CC, KY, HC, CZ, QC1, AH, HP, JW1, JW2, JW3, JW} and references therein). One of the most elegant works is given by He and Xin in \cite{HC, CZ}, and they first realized that the velocity fields played a dominate role in the regularity of the solution to 3D incompressible MHD equations.

\vskip .1in
 The study of the viscous flow in the thin domains started in the seminal paper of  Hale and Raugel \cite{HA}, which was dedicated to a damped hyperbolic equation. More precisely, it is proved that the global attractors are upper semicontinuous. It is shown also that a global attractor exists in the case of the critical sobolev exponent.  In present paper, based on this work about the thin domains mentioned above and motivated by the idea where the  Navier-Stokes equations  converge globally uniformly and strongly to the primitive equations  in \cite{JKL}. We analysis the 3D incompressible MHD equations by the scale technique to derive the primitive equations with magnetic field (PEM) on the thin 3D domains. Because a thin 3D domain is somehow close to a 2D domain, it is natural to use the good properties of the 2D MHD  equations to study the global regularity of strong solutions to  3D PEM  in the thin domains, which is the main idea of our paper.

\vskip .1in
\subsection{Set-up and main results}
 Consider the incompressible three-dimensional  MHD equations in an  $\varepsilon$-dependent thin domain $\Omega_\varepsilon:=M\times(-\varepsilon,\varepsilon)\subset \mathbb{R}^3$, where $\varepsilon>0$ is a very small parameter, and
$M=(0,L_1)\times(0,L_2)$, for two positive constants $L_1$ and $L_2$ of order $O(1)$ with respect to $\varepsilon$. The incompressible three-dimensional anisotropic MHD system is
\begin{equation}\label{1}
 \left\{
\begin{array}{l}
\p_t u + u\cdot \nabla u + \nabla p - \mu\Delta_H u-\nu\partial_z^2 u=b\cdot\nabla b,\\
\p_t b+ u\cdot \nabla b-\kappa\Delta_H b-\sigma\partial_z^2 b=b\cdot\nabla u, \\
\nabla\cdot u =0,  \quad \nabla\cdot b=0,
 \end{array} \right.
 \end{equation}
 where $u=(\U, u_3)$ and $b=(\B, b_3)$ are the velocity field and the magnetic field, respectively. $\U=(u_1, u_2)$ and $\B=(b_1, b_2)$ denote the horizontal velocity field and  magnetic field, respectively, while $u_3$ and $b_3$ stand for the vertical one. The scalar function $p$ is the pressure.  Similar to the case considered in Az\'erad-Guill\'en \cite{AF}, which  emphasized that the anisotropic viscosity hypothesis is fundamental for the derivation of the primitive equations. In our paper, we suppose that  the horizontal   and  the  vertical viscous coefficient $\mu$ and $\nu$ have different orders, that is
  $\mu=O(1)$ and $\nu=O(\varepsilon^2)$. The orders of magnetic diffusivity coefficient $\kappa$ and $\sigma$ are similar to the viscous coefficient. For the sake of simplicity, we set  $\mu=1$ and $\nu=\varepsilon^2,$ similarly, $\kappa=1$ and $\sigma=\varepsilon^2.$ Throughout this paper, we emphasize that the operators $\Delta_H$ act only on horizontal Laplacian, that is $\De_H=\p_x^2+\p_y^2$. Note that it is necessary to consider the above anisotropic viscosity and magnetic diffusivity scaling in the horizontal and vertical directions, so that the MHD equations converge to some specific equations that would be called the primitive equations with magnetic (PEM), as the aspect ratio $\varepsilon$ goes to zero.

\vskip .1in

 We carry out the following scaling transformation to  the equations (\ref{1}) such that the resulting system is defined on a fixed domain independent of $\varepsilon$. To this end, we introduce the new unknowns,
$$u_\varepsilon=(\U_\varepsilon, u_{3,\varepsilon}),\quad b_\varepsilon=(\tilde{b}_\varepsilon, b_{3,\varepsilon}),\quad p_\varepsilon(x,y,z,t)=p(x,y,\varepsilon z,t),$$
$$\U_\varepsilon(x,y,z,t)=\U(x,y,\varepsilon z,t)=(u_1(x,y,\varepsilon z,t), u_2(x,y,\varepsilon z,t)),$$
$$\tilde{b}_\varepsilon(x,y,z,t)=\tilde{b}(x,y,\varepsilon z,t)=(b_1(x,y,\varepsilon z,t), b_2(x,y,\varepsilon z,t)),$$
and
$$u_{3,\varepsilon}(x,y,z,t)=\frac{1}{\varepsilon}u_3(x,y,\varepsilon z,t),\qquad b_{3,\varepsilon}(x,y,z,t)=\frac{1}{\varepsilon}b_3(x,y,\varepsilon z,t).$$

For any $(x,y,z)\in\Omega: = M\times (-1,1) $ and  $t\in (0,\infty)$, then  $u_\varepsilon$, $b_\varepsilon$ and $p_\varepsilon$ satisfy the following  scaled incompressible  MHD equations (SMHD)
\begin{equation}\label{2}
 \left\{
\begin{array}{l}
\p_t \U_\varepsilon + u_\varepsilon\cdot \nabla \U_\varepsilon + \nabla_H p_\varepsilon - b_\varepsilon\cdot\nabla\tilde{b_\varepsilon}-\Delta\U_\varepsilon=0,\\
\varepsilon^2(\partial_t u_{3,\varepsilon}+u_\varepsilon\cdot\nabla u_{3,\varepsilon}-\Delta u_{3,\varepsilon} - b_\varepsilon\cdot\nabla b_{3,\varepsilon})+\partial_z p_\varepsilon=0,\\
\partial_t\tilde{b_\varepsilon}+u_\varepsilon\cdot\nabla\tilde{b_\varepsilon}-\Delta\tilde{b_\varepsilon}-b_{\varepsilon}\cdot\nabla \U_{\varepsilon}=0, \\
\varepsilon^2(\partial_t b_{3,\varepsilon}+u_{\varepsilon}\cdot\nabla b_{3,\varepsilon}-\Delta b_{3,\varepsilon}-b_\varepsilon\cdot\nabla u_{3,\varepsilon})=0,\\
\nabla\C u_\varepsilon=0,\quad \nabla\C b_\varepsilon=0.
 \end{array} \right.
 \end{equation}

 The above equations (\ref{2}) are defined in the fixed domain $\Omega$. Throughout this paper, we set $\na_H$ to denote $(\p_x, \p_y)$. In addition, we equip the system (\ref{2}) with the following initial value conditions and periodic boundary conditions,
 \begin{align}\label{b0}
(\U_\varepsilon,u_{3,\varepsilon})|_{t=0}=(\U_{\varepsilon,0},u_{3,\varepsilon,0}),\qquad (\tilde{b}_\varepsilon,b_{3,\varepsilon})|_{t=0}=(\B_{\varepsilon,0},b_{3,\varepsilon,0}),
\end{align}
and
\begin{align}\label{b1}
\U_\varepsilon,~u_{3,\varepsilon},~\tilde{b}_\varepsilon, ~b_{3,\varepsilon},~ p_\varepsilon  \text{ are periodic in}~x,~y,~z.
\end{align}

It should be noticed that in (\ref{b1}), as well as in all the cases of periodic boundary conditions below, the periods in $x$ and $y$ are $L_1$ and $L_2$, respectively, while that in $z$ is $2$. Furthermore, $(\U_{\varepsilon,0},u_{3,\varepsilon,0})$ and $(\B_{\varepsilon,0},b_{3,\varepsilon,0})$ are given, for simplicity, we suppose in addition that the following symmetry conditions hold
\begin{align}\label{b2}
\U_\varepsilon,~u_{3,\varepsilon}  ~\text{and}~  p_\varepsilon ~\text{are even, odd and odd with respect to}~z,~\text{respectively},
\end{align}
and
\begin{align}\label{b3}
\tilde{b}_\varepsilon,~b_{3,\varepsilon}~\text{are  even and odd with respect to} ~z,~\text{respectively}.
\end{align}

Note that these symmetry conditions are preserved by the dynamics of (SMHD), in other words, they are automatically satisfied as long as they are satisfied initially. So in this article, without further mention, we always assume that the initial horizontal velocity and  magnetic field $\U_0$, $\tilde{b}_0$ satisfy  that
\begin{align*}\label{b4}
\U_0~,\tilde{b}_0~ \text{are  periodic in}~x,~y,~z,~\text{and are even in}~z.
\end{align*}

By taking the limit as  $\varepsilon\rightarrow 0$ in (SMHD) (\ref{2}),  it is natural to obtain the following primitive equations with magnetic field (PEM)
\begin{equation}\label{b5}
 \left\{
\begin{array}{l}
\partial_t \U+u\cdot\nabla \U-\Delta \U-b\cdot\nabla\tilde{b}+\nabla_{H} p = 0,\\
\partial_z p=0,\\
\partial_t\tilde{b}+u\cdot\nabla\tilde{b}-\Delta\tilde{b}-b\cdot\nabla \U=0,\\
\nabla_{H}\cdot \U+\partial_z u_3=0, \qquad \nabla_{H}\cdot\tilde{b}+\partial_z b_3=0.
\end{array} \right.
 \end{equation}

Recalling that we consider the periodic initial-boundary value problem to the (SMHD) equations (\ref{2}), it is clear that one should impose the same boundary conditions and symmetry conditions to the corresponding limiting system (\ref{b5}). However, one only needs to impose the initial conditions on the horizontal velocity field and magnetic field. In fact, since $u_{3,0}$ and $b_{3,0}$ are odd in $z$, we have $u_{3,0}(x,y,0)=b_{3,0}(x,y,0)=0$. Then, $u_{3,0}, b_{3,0}$ can be determined uniquely  by the incompressibility conditions, namely,
\begin{align}\label{3}
u_{3,0}(x,y,z)=-\int_0^z\nabla_H\cdot \U_0(x,y,\xi)\,d\xi,
\end{align}
and
\begin{align}\label{30}
b_{3,0}(x,y,z)=-\int_0^z\nabla_H\C\B_0(x,y,\xi)\,d\xi.
\end{align}

Similarly, $(u_3,b_3)$ can also be  determined uniquely by $(\U, \B)$ via the incompressibility conditions as
\begin{align}\label{t3}
u_{3}(x,y,z,t)=-\int_0^z\nabla_H\cdot \U(x,y,\xi,t)\,d\xi,
\end{align}
and
\begin{align}\label{m1}
b_{3}(x,y,z,t)=-\int_0^z\nabla_H\cdot\tilde{b}(x,y,\xi,t)\,d\xi.
\end{align}

Due to these facts, throughout this paper, concerning the solutions to (\ref{b5}), we only specify the horizontal components $(\U, \B)$, and $(u_3, b_3)$ are determined uniquely  by (\ref{t3}) and (\ref{m1}).


Our main results on  the  global existence and uniqueness (regularity) of strong solutions to the  three-dimensional incompressible PEM (\ref{b5})  without any small assumption on the initial data are stated in the following.

\begin{thm}\label{t1}
 Suppose that $(\U_0,\B_0) \in H^2(\Omega)$,   then there exists a unique  global strong solution $(\U,\B)\in L^{\infty}([0,\infty);H^2(\Om))\cap L^2([0,\infty);H^3(\Om))$ of the  PEM (\ref{b5}), subject to  the boundary and initial conditions (\ref{b0})-(\ref{b3}). Moreover, we have the following estimate,
\begin{align*}
&\sup\limits_{0\leq t< \infty} \|\U\|_{H^2}^2(t)+\sup\limits_{0\leq t<\infty}\|\tilde{b}\|_{H^2}^2(t)+\int_0^\infty(\|\na \U\|_{H^2}^2+\|\partial_t \U\|_{H^1}^2)\,dt\\
&+\int_0^\infty(\|\na\tilde{b}\|_{H^2}^2+\|\partial_t\tilde{b}\|_{H^1}^2)\,dt
\leq C,
\end{align*}
for a constant $C$ depending only on $\|\U_0\|_{H^2}$, $\|\B_0\|_{H^2}$, $L_1$ and $L_2$. Moreover, the unique global strong solution $(\U,\B)$ depends continuously on the initial data.
\end{thm}

\begin{re}
If the initial data  $(\U_0,\B_0)$ belongs to  $H^1(\Om)$, then there exists  a unique global strong  solution to the  PEM (\ref{b5}), it satisfies $(\U,\B)\in L^{\infty}([0,\infty);H^1(\Om))\cap L^2([0,\infty);H^2(\Om))$, which proof  follows directly from the  corollary of Proposition \ref{b9}.
\end{re}

\begin{re}
Generally, if $(\U_0, \B_0)\in H^k$, with $k>2$, there exists a unique global strong solution to the PEM (\ref{b5}), it satisfies $(\U,\B)\in L^\infty([0,\infty), H^k(\Om))\cap L^2([0,\infty);H^{k+1}(\Om))$, then one can show that
\begin{align*}
&\sup\limits_{0\leq t< \infty} \|\U\|_{H^{k}}^2(t)+\sup\limits_{0\leq t<\infty}\|\tilde{b}\|_{H^{k}}^2(t)+\int_0^\infty(\|\na \U\|_{H^{k}}^2+\|\partial_t \U\|_{H^{k-1}}^2)\,dt\\
&+\int_0^\infty(\|\na\tilde{b}\|_{H^{k}}^2+\|\partial_t\tilde{b}\|_{H^{k-1}}^2)\,dt
\leq C,
\end{align*}
for a constant $C$ depending only on $\|\U_0\|_{H^k}$, $\|\B_0\|_{H^k}$, $L_1$ and $L_2$.
\end{re}

\begin{re}
In our forthcoming paper \cite{DLL }, we will rigorously  justify that  the (SMHD) equations (\ref{2})  converge strongly to the PEM (\ref{b5}), globally and uniformly in time.
\end{re}

Here we give the main idea for the proof in this paper. First, in order to obtain the required  uniform $H^2$-norm  estimates, we need to get first-order estimate on $(\U, \B)$. The first-order estimates depend on the $L^\infty([0,\infty),L^2(\Om))$ on $(\p_z\U, \p_z\B)$, therefore, to establish the estimates $(\p_z\U, \p_z\B)$, we will need the estimates of $L^\infty(0,\infty; L^4(\Om))$ on $(\U, \B)$. So the crucial step to prove the global existence of strong solutions is to obtain the estimates of  $L^4$-norm on $(\U,\B)$. One may try to use the standard energy approach to get such  estimates. However, due to the lack of the vertical components estimates $\|u_3\|_4$ and $\|b_3\|_4$ on the left-hand side of the energy inequality, one will encounter two nonlinear terms
  $$\int_\Om b\C\na\B|\U|^2\U\,dxdydz \quad \text{and}\quad \int_\Om b\C\na\U|\B|^2\B\,dxdydz$$
on the right-hand side of the energy inequality which can not be controlled by the $L^4$-norm
$$\sup\limits_{0\leq s\leq t} \|\U\|_4^4(s)+\sup\limits_{0\leq s\leq t} \|\B\|_4^4(s)+2\int_0^t\||\U|\nabla \U\|_2^2+\||\B|\nabla \B\|_2^2\,ds$$
 on the left-hand side. To avoid estimating these two nonlinear terms, we employ the  governing equation for Els\"{a}sser variables. All nonlinear terms in the new equations can be handled via the divergence free conditions.
As we will see in Proposition \ref{dd}, we can successfully achieve the expected  estimates of $L^4$-norm on $(\U,\B)$. As a consequence, based on these estimates, we can obtain other relevant estimates which are sufficient to prove the global existence of strong solution.

\vskip .1in
\subsection{The structure of this paper}

The remainder of this paper is organized as follows.  The Section 2 is dedicated to the basic notations and some Ladyzhenskaya-type inequalities, which will be used in the following sections. In Section 3, the global existence of strong solutions to the three-dimensional incompressible PEM is proved.  In Section 4, we show the continuous dependence on the initial data and the uniqueness of the strong solution.

\section{Preliminaries}
In this section, we introduce the notations used in this paper, and state some Ladyzhenskaya type inequalities  for some kinds of three dimensional integrals, which will be frequently used in the rest of this paper.
\begin{no}

For $q\in[1,\infty]$, we will denote the Lebesgue spaces on the domain $\Omega$ by $L^q=L^q(\Omega)$. For simplicity of notation we will use  $\|\cdot\|_q$  and $\|\cdot\|_{q,M}$ instead of $L^q(\Omega)$ and  $L^q(M)$. For $s\in\mathbb{N}$ the space $H^s(\Omega)$ consists of $f\in L^2(\Omega)$ such that $\nabla^\alpha f\in L^2(\Omega)$ for $|\alpha|\leq s$ endowed with the norm
$$\|f\|_{H^s(\Omega)}=\Big(\sum_{|\alpha|\leq s}\|\na^\alpha f\|_{L^2(\Omega)}^2\Big)^{\frac{1}{2}}.$$
\end{no}
\begin{lem} \label{4}
For convenience, we recall the following Sobolev and Ladyzhenskaya inequalities in $\mathbb{R}^2$, for every $\phi\in H^1(M)$, (see, e.g.,\cite{RAA, PCCF, OAL }),
\begin{align}\label{t24}
\|\phi\|_{L^4(M)}\leq C_0\|\phi\|_{L^2(M)}^{\frac{1}{2}}\|\phi\|_{H^1({M})}^{\frac{1}{2}},
\end{align}
the following Sobolev and Ladyzhenskaya  inequalities in $\mathbb{R}^3$, for every $\phi\in H^1(\Omega)$,
\begin{align}\label{t25}
\|\phi\|_{L^3(\Omega)}\leq C_0\|\phi\|_{L^2(\Omega)}^{\frac{1}{2}}\|\phi\|_{H^1({\Omega})}^{\frac{1}{2}},
\end{align}
and
\begin{align}\label{t26}
\|\phi\|_{L^6(\Omega)}\leq C_0\|\phi\|_{H^1({\Omega})},
\end{align}
we recall the integral version of Minkowsky inequality for the $L^p$ spaces, $p\geq1$. Let $\Omega_1\subset\mathbb{R}^{n_1}$ and $\Omega_2\subset\mathbb{R}^{n_2}$ be two measurable sets, where $n_1$ and $n_2$ are positive integers. Suppose that $f(\xi,\eta)$ is measurable over $\Omega_1\times\Omega_2$. Then,
\begin{align}\label{t27}
\Big[\int_{\Omega_1}\Big(\int_{\Omega_2}|f(\xi,\eta)|\,d\eta\Big)^p\,d\xi\Big]^{\frac{1}{p}}
\leq\int_{\Om_2}\Big(\int_{\Om_1}|f(\xi,\eta)|^p\,d\xi\Big)^{\frac{1}{p}}\,d\eta.
\end{align}
\end{lem}

Next,  we state some Ladyzhenskaya-type inequalities for some kinds of three dimensional integrals.
\begin{lem} \label{b7}
(see \cite[Lemma 2.1]{CCET} ). The following inequalities hold true
\begin{align*}
&\int_M\Big(\int_{-1}^1f(x,y,z)\,dz\Big)\Big(\int_{-1}^1g(x,y,z)h(x,y,z)\,dz\Big)dx dy \\
\leq & C\|f\|_2^{\frac{1}{2}}\Big(\|f\|_2^{\frac{1}{2}}+\|\nabla_{H}f\|_2^{\frac{1}{2}}\Big)\|g\|_2\|h\|_2^{\frac{1}{2}}
\Big(\|h\|_2^{\frac{1}{2}}+\|\nabla_{H}h\|_2^{\frac{1}{2}}\Big),
\end{align*}
and
\begin{align*}
&\int_M\Big(\int_{-1}^1f(x,y,z)\,dz\Big)\Big(\int_{-1}^1g(x,y,z)h(x,y,z)\,dz\Big)dx dy \\
\leq & C\|f\|_2\|g\|_2^{\frac{1}{2}}\Big(\|g\|_2^{\frac{1}{2}}+\|\nabla_H g\|_2^{\frac{1}{2}}\Big)\|h\|_2^{\frac{1}{2}}\Big(\|h\|_2^{\frac{1}{2}}+\|\nabla_H h\|_2^{\frac{1}{2}}\Big),
\end{align*}
for any $f$, $g$, $h$ such that the right-hand sides make sense and are finite, where $C$ is a positive constant depending only on $L_1$ and $L_2$.
\end{lem}
\begin{lem} \label{b8}
(see \cite[Lemma 2.2]{JKL}). Let $\varphi=(\varphi_1,\varphi_2,\varphi_3)$, $\phi$ and $\psi$ be periodic functions with basic domain $\Omega$. Suppose that $\varphi\in H^1(\Omega)$, with $\nabla\cdot\varphi=0$ in $\Omega$, $\int_{\Omega}\varphi dx dy dz =0$, and $\varphi_3|_{z=0}=0$, $\nabla\phi\in H^1(\Omega)$ and $\psi\in L^2(\Omega)$. Denote by $\varphi_H=(\varphi_1,\varphi_2)$ the horizontal components of the function $\varphi$. Then, we have the following estimate
\begin{align*}
\Big|\int_{\Omega}(\varphi\cdot\nabla\phi)\psi\,dxdydz\Big|\leq C\|\nabla\varphi_H\|_2^{\frac{1}{2}}\|\Delta\varphi_H\|_2^{\frac{1}{2}}\|\nabla\phi\|_2^{\frac{1}{2}}\|\Delta\phi\|_2^{\frac{1}{2}}\|\psi\|_2,
\end{align*}
where $C$ is a positive constant depending only on $L_1$, $L_2$.
\end{lem}

\section{A priori estimates on the primitive  equations with magnetic}
\vskip .1in
In this section, we prove \emph{a priori} estimates on the global strong solutions $(\U,u_3)$ and $(\B,b_3)$ to the PEM, we use anisotropic treatments for the PEM to get \emph{a  priori} estimates.
 In Proposition \ref{d1}, the basic energy estimates are shown.
 In Proposition \ref{dd}, we give the $L^\infty(0,\infty;L^4(\Om))$ estimates on $(\U,\B)$.
  Proposition \ref{d2} is devoted to the study of the $L^\infty(0,\infty;L^2(\Omega))$ estimates on ($\partial_z\U$, $\partial_z\tilde{b}$).
  Finally, the $L^\infty(0,\infty;L^2(\Om))$ estimates on  $(\na\U,\na\B)$ and $(\De\U,\De\B)$ are established in Proposition \ref{b9} and Proposition \ref{p20}, respectively.

\begin{prop}\label{d1}
(Basic energy estimates). Suppose that $(\U_0,\tilde{b_0})\in H^1(\Omega)$. Then, we have the following estimates
\begin{align}
\|\U(t)\|_2^2+\|\tilde{b}(t)\|_2^2+2\int_0^t\|\nabla \U\|_2^2\,ds+2\int_0^t\|\nabla\tilde{b}\|_2^2\,ds=\|\U_0\|_2^2+\|\tilde{b}_0\|_2^2,
\end{align}
and
\begin{align}\label{x1}
\|\U(t)\|_2^2+\|\tilde{b}(t)\|_2^2\leq e^{-2C_1t}(\|\U_0\|_2^2+\|\tilde{b}_0\|_2^2).
\end{align}
\end{prop}

\begin{proof}[Proof] Taking the $L^2(\Omega)$ inner product to  the first and the third equation in (\ref{b5}) with $\U$ and $\tilde{b}$, then it follows from integration by parts that
\begin{align*}
\frac{1}{2}\frac{d}{dt}\|\U\|_2^2+\frac{1}{2}\frac{d}{dt}\|\tilde{b}\|_2^2+\|\nabla \U\|_2^2+\|\nabla\tilde{b}\|_2^2=0,
\end{align*}
from which, integrating in $t$ yields the energy identity. In the following estimates, we implicitly use  the  Poincar\'e inequality $\|\U\|_2^2\leq C_1\|\na \U\|_2$, clearly
\begin{align*}
\frac{d}{dt}\|\U\|_2^2+\frac{d}{dt}\|\tilde{b}\|_2^2+2C_1\| \U\|_2^2+2C_1\|\tilde{b}\|_2^2\leq0,
\end{align*}
from which, by the Gronwall inequality, the second conclusion (\ref{x1}) follows.
\end{proof}
Since the high-order estimates depend on the uniform $L^4$-norm estimates of $(\U,\B)$, we first prove these estimates in the following proposition.

\begin{re}
The $L^\infty(0,\infty;L^4(\Omega))$ estimates on $(\U, \B)$ are the foundation of the required $H^2$-norm estimates of Theorem \ref{t1}.
The proof of the  $L^4$-norm estimates on $(\U,\B)$ is not trivial.  A natural starting point is to bound $\|\U\|_4+\|\B\|_4$, via multiplying the first and the third equation of (\ref{b5}) by $|\U|^2\U$ and $|\B|^2\B$, respectively, and integrating the resultant over $\Om$. However, due to the lack of the vertical components $\|u_3\|_4$ and $\|b_3\|_4$ on the left-hand side of the energy inequality, it is very difficult to  bound some of the nonlinear terms directly. More precisely, two of the most troublesome ones are
\begin{align*}
\int_\Om b\C\na\B|\U|^2\U\,dxdydz \quad \text{and} \quad  \int_\Om b\C\na\U|\B|^2\B\,dxdydz,
\end{align*}
 which can not be controlled in terms of $\|\U\|_4+\|\B\|_4$ or the dissipative parts $\||\U|\na\U\|_2^2$ and $\||\B|\na\B\|_2^2$, consequently, we are not able to obtain the uniform bound on the following quantities
\begin{align*}
\sup\limits_{0\leq s\leq t} \|\U\|_4^4(s)+\sup\limits_{0\leq s\leq t} \|\B\|_4^4(s)+2\int_0^t\||\U|\nabla \U\|_2^2+\||\B|\nabla \B\|_2^2\,ds.
\end{align*}

This forces us to avoid estimating these two nonlinear terms. To solve this problem, a crucial step is to introduce the Els\"{a}sser variables $(A=\U+\B, A^*=\U-\B)$ in the equations (\ref{b5}). We add the first equation to the third equation and subtract the third  equation from the first equation in order to obtain a governing equation for Els\"{a}sser variables. Luckily, the new equation (\ref{C1}) serves our purpose perfectly. All nonlinear terms in the equation (\ref{C1}) can be eliminated by the  divergence free conditions   $\na\C u=0$ and $\na\C b=0$. Then, we make use of the $\|\U\|_4+\|\B\|_4 $ can be bounded by $\|\U+\B\|_4$ and $\|\U-\B\|_4$.
\end{re}

\begin{prop}\label{dd}
($L^\infty(0,\infty;L^4(\Omega))$ estimates on $\U$, $\tilde{b}$). Suppose that $(\U_0,\tilde{b_0})\in H^1(\Omega)$. Then, we have the following estimates
\begin{align*}
\sup\limits_{0\leq s\leq t} \|A\|_4^4(s)+\sup\limits_{0\leq s\leq t} \|A^*\|_4^4(s)+\int_0^t\||A|\nabla A\|_2^2\,ds+\int_0^t\||A^*|\nabla A^*\|_2^2\,ds\leq R(0),
\end{align*}
where
$$A=\U+\B,\qquad A^*=\U-\B,$$
and
$$R(0)=exp\{C(\|A_0\|_2^2+\|A_0\|_2^4+\|A_0^*\|_2^2+\|A_0^*\|_2^4)\}(\|A_0^*\|_4^4+\|A_0\|_4^4),$$
for any $t\in [0,\infty)$, where $C$ is a positive constant depending only on $L_1$ and $L_2$.
\end{prop}
\vskip .1in

\begin{proof}[Proof ]
Adding the first equation  and the third equation of (\ref{b5}), we obtain
\begin{align*}
\p_t(\U+\tilde{b})+u\cdot \nabla(\U+\tilde{b})-\Delta(\U+\tilde{b})-b\cdot\nabla(\U+\tilde{b})+\nabla_H p=0,
\end{align*}
subtracting   the third equation from the first equation of (\ref{b5}), yields
\begin{align*}
\p_t(\U-\tilde{b})+u\cdot \nabla(\U-\tilde{b})-\Delta(\U-\tilde{b})+b\cdot\nabla(\U-\tilde{b})+\nabla_H p=0.
\end{align*}

We introduce the Els\"{a}sser variables
$$A=\U+\tilde{b},\qquad A^*=\U-\tilde{b},$$
  and we have the following new formulation for the system (\ref{b5})
\begin{equation}\label{C1}
 \left\{
\begin{array}{l}
\partial_t A+u\cdot\nabla A-\Delta A-b\cdot\nabla{A}+\nabla_{H}p = 0,\\
\partial_t A^*+u\cdot\nabla A^*-\Delta A^*+b\cdot\nabla{A^*}+\nabla_{H}p = 0.\\
\end{array} \right.
 \end{equation}

Multiplying the first equation and the second equation of (\ref{C1}) by $|A|^2A$ and $|A^*|^2A^*$, respectively,  integrating the result over $\Omega$, then it follows from integration by parts that
\begin{align}\label{C2}
&\frac{1}{4}{\ddt}\|A\|_4^4+\frac{1}{4}{\ddt}\|A^*\|_4^4+\int_\Omega|A|^2(|\nabla A|^2+2|\nabla|A||^2)\,dxdydz\nonumber\\
&+\int_\Omega|A^*|^2(|\nabla A^*|^2+2|\nabla|A^*||^2)\,dxdydz\nonumber\\
=&-\int_\Omega|A|^2A\cdot\nabla_Hp\,dxdydz-\int_\Omega|A^*|^2A^*\cdot\nabla_H p\,dxdydz.
\end{align}

 Using Lemma \ref{b7} and the Poincar\'e inequalities, we infer that
\begin{align}\label{c3}
&-\int_\Omega|A|^2A\cdot\nabla_H p(x,y,t)\,dxdydz\nonumber\\
\leq&\int_M\Big(\int_{-1}^1|A^3|\,dz\Big)|\nabla_H p(x,y,t)|\,dxdy\nonumber\\
\leq&C\|\nabla_Hp\|_{2,M}\||A|^2\|_2^{\frac{1}{2}}(\||A|^2\|_2+\|\nabla_H|A|^2\|_2)
^{\frac{1}{2}}\|A\|_2^{\frac{1}{2}}(\|A\|_2+\|\nabla A\|_2)^{\frac{1}{2}}\nonumber\\
\leq&C\|\nabla_H p\|_{2,M}\|A\|_4(\|A\|_4^2+\|\nabla_H|A|^2\|_2)^{\frac{1}{2}}\|A\|_2^{\frac{1}{2}}\|\nabla A\|_2^{\frac{1}{2}}.
\end{align}

Employing the operator $\int_{-1}^1div_H(\cdot)\,dz$ on the   first equation of (\ref{C1}), one obtain
\begin{align}\label{T1}
&\int_{-1}^1\nabla_H\cdot\p_t A\,dz+\int_{-1}^1\nabla_H\cdot(u\cdot\nabla A)\,dz-\int_{-1}^1\nabla_H\cdot(b\cdot\nabla A)\,dz\nonumber\\
&-\int_{-1}^1\nabla_H\cdot(\Delta A)\,dz+\int_{-1}^1\nabla_H\cdot\nabla_H p\,dz=0,
\end{align}
 we denote that
\begin{align*}
&K_1=\int_{-1}^1\nabla_H\cdot\p_t A\,dz,\\
&K_2=\int_{-1}^1\nabla_H\cdot(u\cdot\nabla A)\,dz,\\
&K_3=-\int_{-1}^1\nabla_H\cdot(b\cdot\nabla A)\,dz,\\
&K_4=-\int_{-1}^1\nabla_H\cdot(\Delta A)\,dz,
\end{align*}
and
\begin{align*}
&K_5=\int_{-1}^1\nabla_H\cdot\nabla_H p\,dz.
\end{align*}

Recalling the fact that the periods of $u_3$ and $b_3$  in $z$ are $2$, we obtain
\begin{align*}
K_1=\int_{-1}^1\nabla_H\cdot\p_t A\,dz=-\p_t\int_{-1}^1\p_z u_3+\p_z b_3\,dz=0.
\end{align*}
To bound $K_2$ and $K_3$,  we decompose it into two pieces, respectively
\begin{align*}
K_2+K_3=&\int_{-1}^1\nabla_H\cdot(u\cdot\nabla A)\,dz-\int_{-1}^1\nabla_H\cdot(b\cdot\nabla A)\,dz\\
=&\int_{-1}^1\nabla_H\C(u_3\p_zA)\,dz+\int_{-1}^1\nabla_H\C(\U\cdot\nabla_HA)\,dz\\
&-\int_{-1}^1\nabla_H\cdot(b_3\p_zA)\,dz-\int_{-1}^1\nabla_H\cdot(\B\cdot\nabla_H A)\,dz\\
=&K_{21}+K_{22}+K_{31}+K_{32}.
\end{align*}
To deal with $K_{21}$, integrating by parts and using the divergence free conditions yield that
\begin{align*}
K_{21}=&\int_{-1}^1\nabla_H\cdot(u_3\p_z A)\,dz\\
=&\int_{-1}^1u_3\nabla_H\C\p_z A\,dz+\int_{-1}^1\p_z A\nabla_H u_3\,dz\\
=&\int_{-1}^1u_3\nabla_H\C\p_z(\U+\B)\,dz+\int_{-1}^1\p_z(\U+\B)\nabla_H u_3\,dz\\
=&\int_{-1}^1-u_3\p_z\p_z(u_3+b_3)\,dz+\int_{-1}^1\nabla_H\C(\U+\B)\p_z u_3\,dz\\
=&\int_{-1}^1\p_zu_3\p_z(u_3+b_3)\,dz+\int_{-1}^1\nabla_H\cdot(\U+\B)\p_z u_3\,dz\\
=&\int_{-1}^1\p_zu_3\p_z(u_3+b_3)\,dz-\int_{-1}^1\p_z(u_3+b_3)\p_zu_3\,dz\\
=&0.
\end{align*}

 A similar argument  to  that for $K_{21}$, one gets
\begin{align*}
K_{31}=&-\int_{-1}^1\nabla_H\cdot(b_3\p_z A)\,dz=0.
\end{align*}

On the account of these estimates of $K_2$ and  $K_3$, one has
\begin{align*}
K_2+K_3=&K_{22}+K_{32}\\
=&\int_{-1}^1\nabla_H\cdot(\U\cdot\nabla_HA)\,dz-\int_{-1}^1\nabla_H\cdot(\B\cdot\nabla_H A)\,dz\\
=&\int_{-1}^1\nabla_H\cdot(\U\cdot\nabla_H(\U+\B))\,dz-\int_{-1}^1\nabla_H\cdot(\B\cdot\nabla_H(\U+\B))\,dz\\
=&\int_{-1}^1\nabla_H\cdot(\U\cdot\nabla_H\U)\,dz+\int_{-1}^1\nabla_H\C(\U\cdot\nabla_H\B)\,dz\\
&-\int_{-1}^1\nabla_H\cdot(\B\cdot\nabla_H \U)\,dz-\int_{-1}^1\nabla_H\C(\B\C\na_H\B)\,dz\\
=&\int_{-1}^1\nabla_H\cdot(\nabla_H\cdot(\U\otimes\U))\,dz-\int_{-1}^1\nabla_H\cdot(\nabla_H\cdot(\B\otimes\B))\,dz.
\end{align*}

A similar argument to that for $K_{1}$  yields that
\begin{align*}
K_4=-\int_{-1}^1\nabla_H\cdot(\Delta A)\,dz=0.
\end{align*}

With the aid of the above estimates, the equation (\ref{T1}) becomes as following
\begin{align}\label{7}
\int_{-1}^1\De_H p(x,y,t)\,dz=-\int_{-1}^1\nabla_H\cdot(\nabla_H\cdot(\U\otimes\U))\,dz+\int_{-1}^1\nabla_H\cdot(\nabla_H\cdot(\B\otimes\B))\,dz.
\end{align}

Note that $p$ can be uniquely determined by requiring $\int_\Omega p\,dxdy=0$, and thus, by the elliptic estimates, the pressure can be bounded as following
 \begin{align*}
\|\nabla_Hp\|_{2,M}\leq&\Big\|\int_{-1}^1\nabla_H\cdot(\nabla_H\cdot(\U\otimes\U))\,dz\Big\|_{2,M}+
\Big\|\int_{-1}^1\nabla_H\cdot(\nabla_H\cdot(\B\otimes\B))\,dz\Big\|_{2,M}\\
\leq&C\||\U|\nabla_H\U\|_2+C\||\B|\nabla_H\B\|_2\\
\leq&C\||A|\nabla_HA\|_2.
\end{align*}

Thanks to the above estimate, it follows from (\ref{c3}) and Young inequality that
\begin{align}\label{c4}
&-\int_\Omega|A|^2A\cdot\nabla_H p(x,y,t)\,dxdydz\nonumber\\
\leq& C\||A|\nabla_H A\|_2\|A\|_4(\|A\|_4^2+\|\nabla_H|A|^2\|_2)^{\frac{1}{2}}\|A\|_2^{\frac{1}{2}}\|\nabla A\|_2^{\frac{1}{2}}\nonumber\\
\leq&C(\|A\|_4^2\||A|\nabla_H A\|_2+\|A\|_4\||A|\nabla_H A\|_2^{\frac{3}{2}})\|A\|_2^{\frac{1}{2}}\|\nabla A\|_2^{\frac{1}{2}}\nonumber\\
\leq&\frac{1}{2}\||A|\nabla A\|_2^2+C(\|A\|_2\|\nabla A\|_2+\|A\|_2^2\|\nabla A\|_2^2)\|A\|_4^4.
\end{align}

We deal with the following term  similar to the argument in (\ref{c3}),
\begin{align}\label{C3}
&-\int_\Omega|A^*|^2A^*\cdot\nabla_H p(x,y,t)\,dxdydz\nonumber\\
\leq&\int_M\Big(\int_{-1}^1|A^*|^3\,dz\Big)|\nabla_H p(x,y,t)|\,dxdy\nonumber\\
\leq&C\|\nabla_Hp\|_{2,M}\||A^*|^2\|_2^{\frac{1}{2}}(\||A^*|^2\|_2+\|\nabla_H|A^*|^2\|_2)
^{\frac{1}{2}}\|A^*\|_2^{\frac{1}{2}}(\|A^*\|_2+\|\nabla A^*\|_2)^{\frac{1}{2}}\nonumber\\
\leq&C\|\nabla_H p\|_{2,M}\|A^*\|_4(\|A^*\|_4^2+\|\nabla_H|A^*|^2\|_2)^{\frac{1}{2}}\|A^*\|_2^{\frac{1}{2}}\|\nabla A^*\|_2^{\frac{1}{2}}.
\end{align}

Applying the operator $\int_{-1}^1div_H(\cdot)\,dz$ to the second equation of (\ref{C1}), we obtain
\begin{align}\label{T2}
&\int_{-1}^1\nabla_H\cdot\p_t A^*\,dz+\int_{-1}^1\nabla_H\cdot(u\cdot\nabla A^*)\,dz+\int_{-1}^1\nabla_H\cdot(b\cdot\nabla A^*)\,dz\nonumber\\
&-\int_{-1}^1\nabla_H\cdot(\Delta A^*)\,dz+\int_{-1}^1\nabla_H\cdot\nabla_H p\,dz=0.
\end{align}

 We set that
\begin{align*}
&L_1=\int_{-1}^1\nabla_H\cdot\p_t A^*\,dz,\\
&L_2=\int_{-1}^1\nabla_H\cdot(u\cdot\nabla A^*)\,dz,\\
&L_3=\int_{-1}^1\nabla_H\cdot(b\cdot\nabla A^*)\,dz,\\
&L_4=-\int_{-1}^1\nabla_H\cdot(\Delta A^*)\,dz,
\end{align*}
and
\begin{align*}
L_5=\int_{-1}^1\nabla_H\cdot\nabla_H p\,dz.
\end{align*}

A similar argument to that for $K_1$ and $K_4$, we obtain
\begin{align*}
L_1=\int_{-1}^1\nabla_H\cdot\p_t A^*\,dz=0,
\end{align*}
and
\begin{align*}
L_4=-\int_{-1}^1\nabla_H\cdot(\Delta A^*)\,dz=0.
\end{align*}

The next terms $L_2$ and $L_3$  are split into two parts, respectively
\begin{align*}
L_2+L_3=&\int_{-1}^1\nabla_H\cdot(u\cdot\nabla A^*)\,dz+\int_{-1}^1\nabla_H\cdot(b\cdot\nabla A^*)\,dz\\
=&\int_{-1}^1\nabla_H\C(u_3\p_zA^*)\,dz+\int_{-1}^1\nabla_H\C(\U\cdot\nabla_HA^*)\,dz\\
&\int_{-1}^1\nabla_H\cdot(b_3\p_zA^*)\,dz+\int_{-1}^1\nabla_H\cdot(\B\cdot\nabla_H A^*)\,dz\\
=&L_{21}+L_{22}+L_{31}+L_{32}.
\end{align*}

Along the similar argument to $K_{21}$, yields that
\begin{align*}
L_{21}=&\int_{-1}^1\nabla_H\cdot(u_3\p_z A^*)\,dz=0,
\end{align*}
and
\begin{align*}
L_{31}=&\int_{-1}^1\nabla_H\cdot(b_3\p_z A^*)\,dz=0.
\end{align*}

With the aid of these equality, $L_2$ and $L_3$ can be bounded as follows
\begin{align}\label{5}
L_{2}+L_3=&L_{22}+L_{32}\nonumber\\
=&\int_{-1}^1\nabla_H\cdot(\U\cdot\nabla_HA^*)\,dz+\int_{-1}^1\nabla_H\cdot(\B\cdot\nabla_{H} A^*)\,dz\nonumber\\
=&\int_{-1}^1\nabla_H\cdot(\U\cdot\nabla_H\U)\,dz-\int_{-1}^1\nabla_H\C(\U\cdot\nabla_H\B)\,dz\nonumber\\
&+\int_{-1}^1\nabla_H\cdot(\B\cdot\nabla_H\U)\,dz-\int_{-1}^1\nabla_H\C(\B\cdot\nabla_H\B)\,dz\nonumber\\
=&\int_{-1}^1\nabla_H\cdot(\nabla_H\cdot(\U\otimes\U))\,dz
-\int_{-1}^1\nabla_H\cdot(\nabla_H\cdot(\B\otimes\B))\,dz.
\end{align}

Thanks to the above estimates, the equation (\ref{T2}) becomes as following
\begin{align}\label{8}
\int_{-1}^1\De_H p(x,y,t)\,dz=-\int_{-1}^1\nabla_H\cdot(\nabla_H\cdot(\U\otimes\U))\,dz+\int_{-1}^1\nabla_H\cdot(\nabla_H\cdot(\B\otimes\B))\,dz.
\end{align}

Applying the elliptic estimates, the pressure  can be bounded as following
 \begin{align}\label{6}
\|\nabla_Hp\|_{2,M}\leq&\Big\|\int_{-1}^1\nabla_H\cdot(\nabla_H\cdot(\U\otimes\U))\,dz\Big\|_{2,M}+
\Big\|\int_{-1}^1\nabla_H\cdot(\nabla_H\cdot(\B\otimes\B))\,dz\Big\|_{2,M}\nonumber\\
\leq&C\||\U|\nabla_H\U\|_2+C\||\B|\nabla_H\B\|_2\nonumber\\
\leq&C\||A^*|\nabla_HA^*\|_2.
\end{align}

On the account of (\ref{6}), it follows from (\ref{C3})  that
\begin{align}\label{c5}
&-\int_\Omega|A^*|^2A^*\cdot\nabla_H p(x,y,t)\,dxdydz\nonumber\\
\leq&\frac{1}{2}\||A^*|\nabla A^*\|_2^2+C(\|A^*\|_2\|\nabla A^*\|_2+\|A^*\|_2^2\|\nabla A^*\|_2^2)\|A^*\|_4^4.
\end{align}

Substituting (\ref{c4}) and  (\ref{c5})  into (\ref{C2}), leads to
\begin{align*}
&\frac{d}{dt}(\|A\|_4^4+\|A^*\|_4^4)+2\||A|\na A\|_2^2+2\||A^*|\na A^*\|_2^2\\
\leq &C(\|A\|_2\|\na A\|_2+\|A\|_2^2\|\na A\|_2^2+\|A^*\|_2\|\na A^*\|_2+\|A^*\|_2^2\|\na A^*\|_2^2)(\|A\|_4^4+\|A^*\|_4^4),
\end{align*}
applying the Gronwall inequality to the above inequality, it follows from the H\"{o}lder inequality and Proposition \ref{d1} that
\begin{align*}
&\sup\limits_{0\leq s\leq t} \|A\|_4^4(s)+\sup\limits_{0\leq s\leq t} \|A^*\|_4^4(s)+2\int_0^t\||A|\nabla A\|_2^2\,ds+2\int_0^t\||A^*|\nabla A^*\|_2^2\,ds\\
\leq&exp \Big\{C\Big(\int_0^t(\|A\|_2\|\nabla A\|_2+\|A\|_2^2\|\nabla A\|_2^2)\,ds+\int_0^t(\|A^*\|_2\|\nabla A^*\|_2+\|A^*\|_2^2\|\nabla A^*\|_2^2)\,ds\Big)\Big\}\\
&\times(\|A^*_0\|_4^4+\|A_0\|_4^4)\\
\leq&exp\Big\{ C\Big(\big(\int_0^t\|A\|_2^2\,ds\big)^{\frac{1}{2}}\big(\int_0^t\|\nabla A\|_2^2\,ds\big)^{\frac{1}{2}}+\int_0^t\big(\|A\|_2^2\|\nabla A\|_2^2\big)\,ds\\
&+\big(\int_0^t\|A^*\|_2^2\,ds\big)^{\frac{1}{2}}\big(\int_0^t\|\nabla A^*\|_2^2\,ds\big)^{\frac{1}{2}}+\int_0^t\big(\|A^*\|_2^2\|\nabla A^*\|_2^2\big)\,ds\Big)\Big\}
\times(\|A^*_0\|_4^4+\|A_0\|_4^4)\\
\leq&exp\{C(t^{\frac{1}{2}}e^{-C_1 t}\|A_0\|_2^2+\|A_0\|_2^4+t^{\frac{1}{2}}e^{-C_1 t}\|A^*_0\|_2^2+\|A^*_0\|_2^4)\}(\|A_0^*\|_4^4+\|A_0\|_4^4)\\
\leq&R(0),
\end{align*}
where
\begin{align}
R(0)=exp\{C(\|A_0\|_2^2+\|A_0\|_2^4+\|A_0^*\|_2^2+\|A_0^*\|_2^4)\}(\|A_0^*\|_4^4+\|A_0\|_4^4),
\end{align}
and  $C$ is a positive constant depending only on $L_1$ and $L_2$.

This completes the proof of Proposition \ref{dd} .
\end{proof}
\vskip.1in
 Next, we will give the uniform $L^2$-estimates for $(\p_z\U,\p_z\B)$, which plays a key role in showing the first-order energy estimates on the strong solutions to  the system (\ref{b5}).

\begin{prop}\label{d2}
($L^\infty(0,\infty;L^2(\Omega))$ estimates on $\partial_z\U$, $\partial_z\tilde{b}$). Suppose that $(\U_0,\tilde{b_0})\in H^1(\Omega)$. Then, we have the following estimates
\begin{align*}
&\sup\limits_{0\leq s\leq t} \|\partial_z \U\|_2^2(s)+\sup\limits_{0\leq s\leq t}\|\partial_z \tilde{b}\|_2^2(s)+\int_0^t\|\nabla\partial_z\U\|_2^2+\|\nabla\partial_z\tilde{b}\|_2^2\,ds\\
\leq&K(0),
\end{align*}
where
$$K(0)=\|\partial_z \U_0\|_2^2+\|\partial_z\tilde{b}_0\|_2^2+C(\|\U_0\|_2^2+\|\B_0\|_2^2)(R^2(0)+R(0)).$$
\end{prop}

\begin{proof}[Proof ] Taking the $L^2(\Omega)$ inner product to the first and the third equation of (\ref{b5}) with $-\partial_z^2 \U$ and $-\partial_z^2\tilde{b}$, respectively,  then it follows from integration by parts that
\begin{align}\label{d3}
\frac{1}{2}\frac{d}{dt}\|\partial_z \U\|_2^2+\|\nabla\partial_z \U\|_2^2=&\int_\Omega(\U\cdot\nabla_H \U)\C\partial_z^2 \U+(u_3\partial_z\U)\C\partial_z^2 \U dxdydz\notag\\
&-\int_{\Om}(\tilde{b}\cdot\nabla_H \tilde{b})\C\partial_z^2\U+(b_3\partial_z\tilde{b})\C\partial_z^2 \U\,dxdydz\notag\\
=&I_1+I_2,
\end{align}
and
\begin{align}\label{d4}
\frac{1}{2}\frac{d}{dt}\|\partial_z \tilde{b}\|_2^2+\|\nabla\partial_z \tilde{b}\|_2^2=&\int_\Omega(\U\cdot\nabla_H \tilde{b})\C\partial_z^2 \tilde{b}+(u_3\partial_z\tilde{b})\C\partial_z^2 \tilde{b} \,dxdydz\notag\\
&-\int_{\Om}(\tilde{b}\cdot\nabla_H \U)\C\partial_z^2\tilde{b}+(b_3\partial_z\U)\cdot\partial_z^2\tilde{b}\,dxdydz\notag\\
=&I_3+I_4.
\end{align}

For $I_1$, we obtain
\begin{align*}
I_1=&\int_\Omega(\U\cdot\nabla_H \U)\C\partial_z^2\U+(u_3\partial_z\U)\C\partial_z^2\U\,dxdydz\\
=&-\int_\Omega(\partial_z \U\cdot\nabla_H) \U\cdot\partial_z \U\,dxdydz-\int_\Omega \U\cdot\nabla_H\partial_z\U\cdot\partial_z\U\,dxdydz\\
&-\int_\Omega \partial_z u_3\partial_z\U\cdot\partial_z\U\,dxdydz-\int_\Omega u_3\partial_z^2 \U\cdot\partial_z \U\,dxdydz\\
=&-\int_\Omega(\partial_z \U\cdot\nabla_H) \U\cdot\partial_z \U\,dxdydz+\frac{1}{2}\int_\Omega\nabla_H\cdot \U(\partial_z \U)^2\,dxdydz,\\
&+\int_\Omega\nabla_H\cdot \U\partial_z \U\C\partial_z \U\,dxdydz+\frac{1}{2}\int_\Omega\partial_z u_3(\partial_z \U)^2\,dxdydz\\
=&\int_\Omega[(-\partial_z \U\cdot\nabla_H) \U+\nabla_H\cdot \U\partial_z \U]\C\partial_z \U\,dxdydz,
\end{align*}
where the last equality follows from the fact that  $\na_{H}\C\U+\p_z u_3=0$.

The term $I_2$ can be split into four parts
\begin{align*}
I_{2}=&-\int_\Omega(\tilde{b}\cdot\nabla_H\tilde{b})\C\partial_z^2 \U\,dxdydz-\int_\Omega(b_3\partial_z\tilde{b})\C\partial_z^2\U\,dxdydz\\
=&\int_\Omega (\partial_z\tilde{b}\cdot\nabla_H)\tilde{b}\C\partial_z \U\,dxdydz+\int_\Omega\tilde{b}\cdot\partial_z\nabla_H\tilde{b}\C\partial_z \U\,dxdydz\\
+&\int_{\Om}\p_z b_3\p_z\B\p_z\U\,dxdydz+\int_{\Om}b_3\p_z^2\B\p_z\U\,dxdydz\\
=&I_{21}+I_{22}+I_{23}+I_{24}.
\end{align*}
A similar argument to that for $I_1$,  applying integration by parts and the incompressible conditions give that
\begin{align*}
I_{3}=&\int_\Omega(\U\cdot\nabla_H\tilde{b})\C\partial_z^2\tilde{b}\,dxdydz+\int_\Omega(u_3\p_z\B)
\C\p_z^2\B\,dxdydz\\
=&-\int_\Omega[(\partial_z \U\C\nabla_H)\tilde{b}-\nabla_H\cdot \U\partial_z \tilde{b}]\C\p_z\tilde{b}\,dxdydz.
\end{align*}

Similarly, $I_4$ can  break it down, namely,
\begin{align*}
I_{4}=&-\int_\Omega(\tilde{b}\cdot\nabla_H \U)\C\partial_z^2 \tilde{b}\,dxdydz-\int_\Omega(b_3\partial_z \U)\C\partial_z^2\tilde{b}\,dxdydz\\
=&\int_\Omega (\partial_z\tilde{b}\cdot\nabla_H) \U\partial_z \tilde{b}\,dxdydz+\int_\Omega\tilde{b}\cdot\partial_z\nabla_H \U\C\partial_z \tilde{b}\,dxdydz\\
&+\int_\Omega \partial_zb_3\partial_z\U\partial_z \tilde{b}\,dxdydz+\int_\Omega b_3\partial_z^2\U\partial_z \tilde{b}\,dxdydz\\
=&I_{41}+I_{42}+I_{43}+I_{44}.
\end{align*}

For the sake of simplicity, we sum up the following four terms, then
\begin{align}\label{9}
&I_{24}+I_{44}+I_{23}+I_{43}\nonumber\\
=&2\int_\Omega\partial_z b_3\partial_z\tilde{b}\partial_z\U\,dxdydz+\int_\Omega b_3\partial_z^2\tilde{b}\partial_z \U\,dxdydz+\int_\Om b_3\p_z^2\U\p_z\B\,dxdydz\nonumber\\
=&2\int_\Om\p_zb_3\p_z\B\p_z\U\,dxdydz-\int_\Om\p_zb_3\p_z\B\p_z\U\,dxdydz
-\int_\Om b_3\p_z\B\p_z^2\U\,dxdysz
+\int_\Om b_3\p_z^2\U\p_z\B\,dxdydz\nonumber\\
=&2\int_\Om\p_zb_3\p_z\B\p_z\U\,dxdydz-\int_\Om\p_zb_3\p_z\B\p_z\U\,dxdydz\nonumber\\
=&-\int_\Om\na_H\C\B\p_z\B\p_z\U\,dxdydz,
\end{align}
which implies that
\begin{align*}
I_{2}+I_{4}=&\int_\Omega[(\p_z\B\cdot\nabla_H)\B-\nabla_H\cdot\B\p_z\B]\cdot\p_z \U\,dxdydz+\int_\Omega\tilde{b}\cdot\partial_z\nabla_H\tilde{b}\cdot\partial_z \U\,dxdydz\\
&+\int_\Omega\partial_z\tilde{b}\cdot\nabla_H \U\cdot\partial_z\tilde{b}\,dxdydz
+\int_\Omega\tilde{b}\cdot\partial_z\nabla_H \U\cdot\partial_z\tilde{b}\,dxdydz.
\end{align*}

Substituting these  estimates of $I_1$-$I_4$ into  (\ref{d3}) and (\ref{d4}) leads to
\begin{align}\label{d5}
&\frac{1}{2}\frac{d}{dt}\|\partial_z \U\|_2^2+\|\nabla\partial_z \U\|_2^2+\frac{1}{2}\frac{d}{dt}\|\partial_z\tilde{b}\|_2^2+\|\nabla\partial_z\tilde{b}\|_2^2\notag\\
=&-\int_\Omega[(\partial_z \U\cdot\nabla_H)\U-\nabla_H\cdot \U\partial_z \U]\cdot\partial_z \U dxdydz-\int_\Omega[(\p_z\U\cdot\nabla_H)\B-\nabla_H\cdot\U\p_z\B]\C\p_z\B\,dxdydz\notag\\
&+\int_\Omega[(\p_z\B\cdot\nabla_H)\B-\nabla_H\cdot\B\p_z\B]\C\p_z\U\,dxdydz
+\int_\Omega\tilde{b}\cdot\partial_z\nabla_H\tilde{b}\cdot\partial_z \U\,dxdydz\notag\\
&+\int_\Omega\partial_z\tilde{b}\cdot\nabla_H \U\cdot\partial_z\tilde{b}\,dxdydz
+\int_\Omega\tilde{b}\cdot\partial_z\nabla_H \U\cdot\partial_z\tilde{b}\,dxdydz\notag\\
\leq&J_1+J_2+J_3+J_4+J_5+J_6,
\end{align}
and it follows from integration by parts that
\begin{align*}
|J_{1}|=&\int_\Omega[(\p_z \U\cdot\na_H)\U-\na_H\C\U\p_z\U]\C\p_z \U\,dxdydz\\
=&\int_\Omega\p_z\U_i\p_i\U_j\p_z\U_j\,dxdydz-\int_\Om\p_i\U_i\p_z\U_j\p_z\U_j\,dxdydz\\
=&-\int_\Om\U_i\p_i\p_z\U_j\p_z\U_j\,dxdydz-\int_\Om\U_i\p_i\U_j\p_z\p_z\U_j\,dxdydz-\int_\Om\p_i\U_i
\p_z\U_j\p_z\U_j\,dxdydz\notag\\
=&\int_\Om\p_i\U_i\p_z\U_j\p_z\U_j\,dxdydz+\int_\Om\U_i\p_z\U_j\p_i\p_z\U_j\,dxdydz\\
&-\int_\Om\U_i\p_i\U_j\p_z\p_z\U_j\,dxdydz-\int_\Om\p_i\U_i\p_z\U_j\p_z\U_j\,dxdydz\\
=&\int_\Om\U_i\p_z\U_j\p_i\p_z\U_j\,dxdydz+\int_\Om\p_z\U_i\p_i\U_j\p_z\U_j\,dxdydz+
\int_\Om\U_i\p_z\p_i\U_j\p_z\U_j\,dxdydz\notag\\
=&\int_\Om\U_i\p_z\U_j\p_i\p_z\U_j\,dxdydz-\int_\Om\p_i\p_z\U_i\U_j\p_z\U_j\,dxdydz\\
&-\int_\Om\p_z\U_i\U_j\p_i\p_z\U_j\,dxdydz+\int_\Om\U_i\p_z\p_i\U_j\p_z\U_j\,dxdydz\\
\leq&4\int_\Om|\p_z\U||\na_H\p_z\U||\U|\,dxdydz,
\end{align*}
employing the H\"{o}lder and Young inequalities, we infer that
\begin{align}\label{p2}
|J_1|\leq&4\int_\Om|\p_z\U||\na_H\p_z\U|\|\U|\,dxdydz\notag\\
\leq&C\|\p_z\U\|_4\|\U\|_4\|\na_H\p_z\U\|_2\notag\\
\leq&C\|\U\|_4\|\p_z\U\|_2^{\frac{1}{4}}\|\na\p_z\U\|_2^{\frac{7}{4}}\notag\\
\leq&\frac{1}{14}\|\na\p_z\U\|_2^2+C\|\U\|_4^8\|\p_z\U\|_2^2.
\end{align}

It follows from integration by parts that
\begin{align*}
|J_{2}|=&\int_\Omega[(\p_z \U\cdot\na_H)\B-\na_H\C\U\p_z\B]\C\p_z \B\,dxdydz\\
\leq&2\int_\Om|\U||\p_z\B||\na_H\p_z\B|\,dxdydz+\int_\Om|\na_H\p_z\U||\B||\p_z\B|\,dxdydz
+\int_\Om|\p_z\U||\B||\p_z\na_H\B|\,dxdydz\\
\leq &J_{21}+J_{22}+J_{23}.
\end{align*}

Firstly, the estimate for $J_{21}$ is given as follows. By the H\"{o}lder and Young  inequalities, we deduce
\begin{align}\label{p3}
J_{21}=&2\int_\Om|\U||\p_z\B||\na_H\p_z\B|\,dxdydz\notag\\
\leq&C\|\p_z\B\|_4\|\U\|_4\|\na_H\p_z\B\|_2\notag\\
\leq&C\|\U\|_4\|\p_z\B\|_2^{\frac{1}{4}}\|\na_H\p_z\B\|_2^{\frac{7}{4}}\notag\\
\leq&\frac{1}{16}\|\na\p_z\B\|_2^2+C\|\U\|_4^8\|\p_z\B\|_2^2.
\end{align}

Secondly, we can  estimate $J_{22}$ as follows
\begin{align}\label{p4}
J_{22}=&\int_\Om|\na_H\p_z\U||\B||\p_z\B|\,dxdydz\notag\\
\leq &C\|\na_H\p_z\U\|_2\|\B\|_4\|\p_z\B\|_4\notag\\
\leq&C\|\na\p_z\U\|_2\|\B\|_4\|\p_z\B\|_2^{\frac{1}{2}}\|\na\p_z\B\|_2^{\frac{1}{2}}\notag\\
\leq&\frac{1}{14}\|\na\p_z\U\|_2^2+\frac{1}{16}\|\na\p_z\B\|_2^2+C\|\B\|_4^4\|\p_z\B\|_2^2,
\end{align}
 and applying the H\"{o}lder and Young  inequalities once again
\begin{align}\label{p5}
J_{23}=&\int_\Om|\p_z\U||\B||\p_z\na_H\B|\,dxdydz\notag\\
\leq &C\|\B\|_4\|\p_z\U\|_4\|\na_H\p_z\B\|_2\notag\\
\leq&C\|\B\|_4\|\p_z\U\|_2^{\frac{1}{2}}\|\na_H\p_z\U\|_2^{\frac{1}{2}}\|\na_H\p_z\B\|_2\notag\\
\leq&\frac{1}{14}\|\na\p_z\U\|_2^2+\frac{1}{16}\|\na\p_z\B\|_2^2+C\|\B\|_4^4\|\p_z \U\|_2^2.
\end{align}

For the estimate of $J_3$, we  break it down
\begin{align*}\label{p6}
J_{3}=&\int_\Om[(\p_z\B\C\na_H)\B-\na_H\C\B\p_z\B]\C\p_z\U\,dxdydz\\
\leq&\int_\Om\p_z\B_i\p_i\B_j\p_z\U_j\,dxdydz-\int_\Om\p_i\B_i\p_z\B_j\p_z\U_j\,dxdydz\\
\leq&-\int_\Om\B_i\p_z\B_j\p_i\p_z\U_j\,dxdydz-\int_\Om\p_i\p_z\B_i\B_j\p_z\U_j\,dxdydz\\
&-\int_\Om\p_z\B_i\B_j\p_i\p_z\U_j\,dxdydz+\int_\Om\B_i\p_i\p_z\B_j\p_z\U_j\,dxdydz\\
\leq&2\int_\Om|\na_H\p_z\U||\B||\p_z\B|\,dxdydz+2\int_\Om|\B||\na_H\p_z\B||\p_z\U|\,dxdydz\\
\leq&J_{31}+J_{32}.
\end{align*}

To bound $J_{31}$,   thanks to  H\"{o}lder and Young  inequalities, one has
\begin{align}
J_{31}=&2\int_\Om|\na_H\p_z\U||\B||\p_z\B|\,dxdydz\notag\\
\leq&2\|\B\|_4\|\p_z\B\|_4\|\na_H\p_z\U\|_2\notag\\
\leq&C\|\B\|_4\|\p_z\B\|_2^{\frac{1}{2}}\|\na_H\p_z\B\|_2^{\frac{1}{2}}\|\na_H\p_z\U\|_2\notag\\
\leq&\frac{1}{14}\|\na\p_z\U\|_2^2+\frac{1}{16}\|\na\p_z\B\|_2^2+C\|\B\|_4^4\|\p_z\B\|_2^2,
\end{align}
for $J_{32}$, we  apply the H\"{o}lder and Young  inequalities once again
\begin{align}\label{p8}
J_{32}=&2\int_\Om|\B||\na_H\p_z\B||\p_z\U|\,dxdydz\notag\\
\leq&C\|\B\|_4\|\p_z\U\|_4\|\na_H\p_z\B\|_2\notag\\
\leq&C\|\B\|_4\|\p_z\U\|_2^{\frac{1}{2}}\|\na_H\p_z\U\|_2^{\frac{1}{2}}\|\na_H\p_z\B\|_2\notag\\
\leq&\frac{1}{14}\|\na\p_z\U\|_2^2+\frac{1}{16}\|\na\p_z\B\|_2^2+C\|\B\|_4^4\|\p_z\U\|_2^2.
\end{align}

A similar argument to that for $J_{32}$ yields
\begin{align}\label{p10}
J_{4}=&\int_\Om|\B||\p_z\U||\na_H\p_z\B|\,dxdydz\notag\\
\leq&C\|\B\|_4\|\p_z\U\|_4\|\na_H\p_z\B\|_2\notag\\
\leq&C\|\B\|_4\|\p_z\U\|_2^{\frac{1}{2}}\|\na_H\p_z\U\|_2^{\frac{1}{2}}\|\na_H\p_z\B\|_2\notag\\
\leq&\frac{1}{14}\|\na\p_z\U\|_2^2+\frac{1}{16}\|\na\p_z\B\|_2^2+C\|\B\|_4^4\|\p_z\U\|_2^2,
\end{align}
and
\begin{align}\label{p11}
J_{5}=&\int_\Om\p_z\B\C\na_H\U\p_z\B\,dxdydz\notag\\
\leq&2\int_\Om|\na_H\p_z\B||\U||\p_z\B|\,dxdydz\notag\\
\leq&C\|\U\|_4\|\p_z\B\|_4\|\na_H\p_z\B\|_2\notag\\
\leq&C\|\U\|_4\|\p_z\B\|_2^{\frac{1}{4}}\|\na\p_z\B\|_2^{\frac{7}{4}}\notag\\
\leq&\frac{1}{16}\|\na\p_z\B\|_2^2+C\|\U\|_4^8\|\p_z\B\|_2^2.
\end{align}

For the last term $J_6$, we get that
\begin{align}\label{p12}
J_{6}=&\int_\Om|\B||\p_z\na_H\U||\p_z\B|\,dxdydz\notag\\
\leq&C\|\B\|_4\|\p_z\B\|_4\|\p_z\na_H\U\|_2\notag\\
\leq&C\|\B\|_4\|\p_z\B\|_2^{\frac{1}{2}}\|\na\p_z\B\|_2^{\frac{1}{2}}\|\p_z\na_H\U\|_2\notag\\
\leq&\frac{1}{14}\|\na\p_z\U\|_2^2+\frac{1}{16}\|\na\p_z\B\|_2^2+C\|\B\|_4^4\|\p_z\B\|_2^2.
\end{align}

With the aid of the above  estimates, substituting (\ref{p2})-(\ref{p12}) into (\ref{d5}), yields
\begin{align}\label{p13}
&\frac{1}{2}\ddt(\|\p_z\U\|_2^2+\|\p_z\B\|_2^2)+\frac{1}{2}(\|\na\p_z\U\|_2^2+\|\na\p_z\B\|_2^2)\notag\\
\leq &C(\|\U\|_4^8+\|\U\|_4^4+\|\B\|_4^4)(\|\p_z\U\|_2^2+\p_z\B\|_2^2),
\end{align}
from which, by the Gronwall inequality, it follows from  Propositions \ref{d1} and  Propositions \ref{dd}, we have
\begin{align}\label{p14}
&\sup\limits_{0\leq s\leq t} \|\p_z\U\|_4^4(s)+\sup\limits_{0\leq s\leq t} \|\p_z\B\|_4^4(s)+\int_0^t\|\na\p_z\U\|_2^2\,ds+\int_0^t\|\na\p_z\B\|_2^2\,ds\notag\\
\leq&\|\p_z\U_0\|_2^2+\|\p_z\B_0\|_2^2+C(\|\U\|_4^8+\|\U\|_4^4+\|\B\|_4^4)\int_0^t\|
\p_z\U\|_2^2+\|\p_z\B\|_2^2\,ds\notag\\
\leq&K(0),
\end{align}
 where
 $$K(0)=\|\p_z\U_0\|_2^2+\|\p_z\B_0\|_2^2+C(\|\U_0\|_2^2+\|\B_0\|_2^2)(R^2(0)+R(0)).$$
 This  completes the proof of the $L^\infty(0,\infty;L^2(\Omega))$ estimates on $\partial_z\U$  and $\partial_z\tilde{b}$.
\end{proof}
\vskip .1in
In this subsection, we work on the first-order energy estimates to the strong solutions of  the system (\ref{b5}), subject to the boundary and initial conditions (\ref{b0})-(\ref{b3}). In particular, we have shown that the growth of the $H^1$-norms of $(\U,\B)$ is not faster than a uniform constant, which depend only on $\|\U_0\|_{H^1}$, $\|\B_0\|_{H^1}$, $L_1$, $L_2$.

\begin{prop}\label{b9}
(First-order energy estimates). Suppose that $(\U_0,\tilde{b_0})\in H^1(\Omega)$, then, we have the following estimates
\begin{align}
&\sup\limits_{0\leq s\leq t} \|\nabla \U\|_2^2(s)+\sup\limits_{0\leq s\leq t}\|\nabla\tilde{b}\|_2^2(s)+\frac{1}{2}\int_0^t(\|\Delta \U\|_2^2+\|\partial_t \U\|_2^2)ds+\frac{1}{2}\int_0^t(\|\Delta\tilde{b}\|_2^2+\|\partial_t\tilde{b}\|_2^2)ds\notag\\
\leq&H(0),
\end{align}
where
$$H(0)=(\|\na\U_0\|_2^2+\|\na\B_0\|_2^2)exp\{C(\|\U_0\|_2^4+\|\B_0\|_2^4+K^2(0))\}.$$
\end{prop}
\vskip .1in
\begin{proof}[Proof] Taking the $L^2(\Omega)$ inner product to the first and the third equation of (\ref{b5}) with $\partial_t \U-\Delta \U$ and $\partial_t\tilde{b}-\Delta\tilde{b}$, respectively,  then it follows from integration by parts that
\begin{align}\label{50}
&\frac{d}{dt}\|\nabla \U\|_2^2+\|\Delta \U\|_2^2+\|\partial_t \U\|_2^2+\frac{d}{dt}\|\nabla\tilde{b}\|_2^2+\|\Delta\tilde{b}\|_2^2+\|\partial_t\tilde{b}\|_2^2 \notag\\
=&\int_{\Omega}[(\U\cdot\nabla_H)\U+u_3\partial_z \U]\cdot(\Delta \U-\partial_t \U)\,dxdydz+\int_{\Om}[(\U\cdot\nabla_H)\tilde{b}+u_3\partial_z\tilde{b}]\C(\Delta\tilde{b}-\partial_t\tilde{b})\,dxdydz\notag\\
&+\int_{\Om}[(\tilde{b}\cdot\nabla_H)\tilde{b}+b_3\partial_z \tilde{b}]\cdot(\partial_t \U-\Delta \U)\,dxdydz+\int_{\Om}[(\tilde{b}\cdot\nabla_H)\U+b_3\partial_z\U]\C(\partial_t\tilde{b}-\Delta\tilde{b})\,dxdydz\notag\\
=&M_1+M_2+M_3+M_4.
\end{align}

To estimate $M_1$, we split it into two terms
\begin{align*}
M_{11}=&\int_\Om(\U\C\na_H)\U\C(\De\U-\p_t\U)\,dxdydz+\int_\Om u_3\p_z\U\C(\De\U-\p_t\U)\,dxdydz\\
=&M_{11}+M_{12}.
\end{align*}

Applying  Lemma \ref{b7}, it follows from  the  Poincar\'e and Young inequalities that
\begin{align} \label{a1}
M_{11}\leq &\int_\Omega(\U\cdot\nabla_H)\U\cdot(\Delta \U-\partial_t \U)\,dxdydz \notag\\
\leq&\int_M\Big(\int_{-1}^1(|\U|+|\partial_z\U|)dz\Big)\Big(\int_{-1}^1|\nabla_H \U|(|\Delta \U|+|\partial_t \U|)dz\Big)\,dxdy \notag\\
\leq&C\Big[\|\U\|_2^{\frac{1}{2}}(\|\U\|_2+\|\nabla_H \U\|_2)^{\frac{1}{2}}+\|\partial_z \U\|_2^{\frac{1}{2}}(\|\partial_z\U\|_2+\|\nabla_H\partial_z \U\|_2)^{\frac{1}{2}} \notag\Big]\\
&\times\|\nabla_H \U\|_2^{\frac{1}{2}}(\|\nabla_H \U\|_2+\|\nabla_H^2 \U\|_2)^{\frac{1}{2}}(\|\Delta \U\|_2+\|\partial_t \U\|_2) \notag\\
\leq& C(\|\U\|_2^{\frac{1}{2}}\|\nabla \U\|_2^{\frac{1}{2}}+\|\partial_z\U\|_2^{\frac{1}{2}}\|\nabla\partial_z \U\|_2^{\frac{1}{2}})\|\nabla_H \U\|_2^{\frac{1}{2}}\|\Delta \U\|_2^{\frac{1}{2}}(\|\Delta \U\|+\|\partial_t \U\|_2) \notag\\
\leq&\frac{1}{12}\|\Delta \U\|_2^2+\frac{1}{8}\|\partial_t \U\|_2^2+C(\|\U\|_2^2\|\nabla \U\|_2^2+\|\partial_z \U\|_2^2\|\nabla\partial_z \U\|_2^2)\|\nabla \U\|_2^2.
\end{align}

Since $u_3$ is odd in $z$, it has $u_3|_{z=0}=0$, and thus
\begin{align*}
u_3(x,y,z,t)=\int_0^z \partial_z u_3(x,y,\xi,t)d\xi=-\int_0^z \nabla_H\cdot \U(x,y,\xi,t)d\xi.
\end{align*}

Thanks to this fact, it follows from Lemma \ref{b7} that
\begin{align}\label{a2}
M_{12}\leq&\int_M u_3\partial_z \U\C(\Delta \U-\partial_t \U)\,dx dy dz\notag\\
\leq&\int_\Omega\Big(\int_{-1}^1|\nabla_H\cdot \U|dz\Big)\Big(\int_{-1}^1|\partial_z \U|(|\Delta \U|+|\partial_t \U|)dz\Big)\,dx dy\notag\\
\leq& C \|\nabla_H \U\|_2^{\frac{1}{2}}(\|\nabla_H \U\|_2+\|\nabla_H^2\U\|_2)^{\frac{1}{2}}\|\partial_z \U\|_2^{\frac{1}{2}}\notag\\
&\times(\|\partial_z \U\|_2+\|\nabla_H \partial_z \U\|_2)^{\frac{1}{2}}(\|\Delta \U\|_2+\|\partial_t \U\|_2)\notag\\
\leq& C\|\nabla_H \U\|_2^{\frac{1}{2}}\|\Delta_H \U\|_2^{\frac{1}{2}}\|\partial_z \U\|_2^{\frac{1}{2}}\|\nabla\partial_z \U\|_2^{\frac{1}{2}}(\|\Delta \U\|_2+\|\partial_t \U\|_2)\notag\\
\leq&\frac{1}{12}\|\Delta \U\|_2^2+\frac{1}{8}\|\partial_t \U\|_2^2+C\|\partial_z \U\|_2^2\|\nabla\partial_z \U\|_2^2\|\nabla \U\|_2^2.
\end{align}

For $M_2$, we decompose it into two pieces,
\begin{align*}
M_2=&\int_\Omega[(\U\cdot\nabla_H)\tilde{b}+u_3\partial_z\tilde{b}]\C(\Delta\tilde{b}-\partial_t\tilde{b})dxdydz\\
=&\int_{\Omega} (\U\cdot\nabla_H) \tilde{b}\cdot(\Delta\tilde{b}-\partial_t\tilde{b})\,dxdydz +\int_\Omega u_3\partial_z\tilde{b}\C(\Delta\tilde{b}-\partial_t\tilde{b})\,dx dy dz\\
=&M_{21}+M_{22}.
\end{align*}

By the  Lemma \ref{b7} and the Poincar\'e and Young inequalities, we deduce
\begin{align} \label{a3}
M_{21}=&\int_{\Omega} (\U\cdot\nabla_H) \tilde{b}\cdot(\Delta\tilde{b}-\partial_t\tilde{b})\,dxdydz \notag\\
\leq&\int_M\Big(\int_{-1}^1(|\U|+|\partial_z \U|)dz\Big)\Big(\int_{-1}^1|\nabla_H \tilde{b}|(|\Delta\tilde{b}|+|\partial_t\tilde{b}|)dz\Big)\,dxdy \notag\\
\leq& C\Big[\|\U\|_2^{\frac{1}{2}}(\|\U\|_2+\|\nabla_H \U\|_2)^{\frac{1}{2}}+\|\partial_z \U\|_2^{\frac{1}{2}}(\|\partial_z \U\|_2+\|\nabla_H \partial_z \U\|_2)^{\frac{1}{2}}\Big] \notag\\
&\times\|\nabla_H\tilde{b} \|_2^{\frac{1}{2}}(\|\nabla_H\tilde{b} \|_2+\|\nabla_H^2\tilde{b} \|_2)^{\frac{1}{2}}(\|\Delta\tilde{b}\|_2+\|\partial_t\tilde{b}\|_2) \notag\\
\leq& C\Big(\|\U\|_2^{\frac{1}{2}}\|\nabla \U\|_2^{\frac{1}{2}}+\|\partial_z \U\|_2^{\frac{1}{2}}\|\nabla\partial_z \U\|_2^{\frac{1}{2}}\Big)\|\nabla_H\tilde{b}\|_2^{\frac{1}{2}}\|\Delta\tilde{b}\|_2^{\frac{1}{2}}(\|\Delta\tilde{b}\|_2+\|\partial_t \tilde{b}\|_2) \notag\\
\leq&\frac{1}{12}\|\Delta\tilde{b}\|_2^2+\frac{1}{8}\|\partial_t \tilde{b}\|_2^2+C(\|\U\|_2^2\|\nabla \U\|_2^2+\|\partial_z \U\|_2^2\|\nabla\partial_z \U\|_2^2)\|\nabla\tilde{b}\|_2^2.
\end{align}

Similar to the argument in (\ref{a2}), it follows from Lemma \ref{b7}, the Poincar\'e and Young inequalities that
\begin{align}\label{a4}
M_{22}=&\int_\Omega u_3\partial_z\tilde{b}\C(\Delta\tilde{b}-\partial_t\tilde{b})\,dx dy dz\notag\\
\leq&\int_M\Big(\int_{-1}^1|\nabla_H\cdot \U|dz\Big)\Big(\int_{-1}^1|\partial_z\tilde{b}|(|\Delta\tilde{b}|+|\partial_t\tilde{b}|)dz\Big)\,dxdy\notag\\
\leq&\|\nabla_H \U\|_2^{\frac{1}{2}}(\|\nabla_H \U\|_2+\|\nabla_H^2 \U\|_2)^{\frac{1}{2}}\|\partial_z\tilde{b}\|_2^{\frac{1}{2}}\notag\\
&\times(\|\partial_z\tilde{b}\|_2+\|\nabla_H \partial_z \tilde{b}\|_2)^{\frac{1}{2}}(\|\Delta\tilde{b}\|_2+\|\partial_t\tilde{b}\|_2)\notag\\
\leq&\|\nabla_H \U\|_2^{\frac{1}{2}}\|\Delta_H \U\|_2^{\frac{1}{2}}\|\partial_z\tilde{b}\|_2^{\frac{1}{2}}\|\nabla\partial_z\tilde{b}\|_2^{\frac{1}{2}}(\|\Delta\tilde{b}\|_2+\|\partial_t\tilde{b}\|_2)\notag\\
\leq&\frac{1}{12}\|\Delta\tilde{b}\|_2^2+\frac{1}{8}\|\partial_t\tilde{b}\|_2^2+\frac{1}{12}\|\Delta\tilde{u}\|_2^2+C\|\partial_z\tilde{b}\|_2\|\nabla\partial_z\tilde{b}\|_2^2\|\nabla \U\|_2.
\end{align}

To deal with $M_3$, we break it down
\begin{align*}
M_3=&\int_\Omega[(\tilde{b}\cdot\nabla_H)\tilde{b}+b_3\partial_z \tilde{b}]\C(\partial_t \U-\Delta \U)\,dxdydz\\
=&M_{31}+M_{32},
\end{align*}
recalling the Lemma \ref{b7}, the Poincar\'e and Young inequalities,  one can easily check
\begin{align}\label{a5}
M_{31}=&\int_{\Omega} (\tilde{b}\cdot\nabla_H) \tilde{b}\cdot(\Delta \U-\partial_t \U)
\,dxdydz\notag\\
\leq&\int_M\Big(\int_{-1}^1(|\tilde{b}|+|\partial_z\tilde{b}|)dz\Big)\Big(\int_{-1}^1|\nabla_H \tilde{b}|(|\Delta \U|+|\partial_t \U|)\,dz\Big)\,dxdy\notag\\
\leq &C\Big[\|\tilde{b}\|_2^{\frac{1}{2}}(\|\tilde{b}\|_2+\|\nabla_H\tilde{b}\|_2)^{\frac{1}{2}}+
\|\p_z\tilde{b}\|_2^{\frac{1}{2}}(\|\partial_z\tilde{b}\|_2+\|\nabla_H\partial_z\tilde
{b}\|_2)^{\frac{1}{2}}\Big]\notag\\
&\times\|\nabla_H\tilde{b}\|_2^{\frac{1}{2}}(\|\nabla_H\tilde{b}\|_2+\|\nabla_H^2\tilde{b}\|_2)
^{\frac{1}{2}}(\|\Delta \U\|_2+\|\partial_t \U\|_2)\notag\\
\leq &C\Big(\|\tilde{b}\|_2^{\frac{1}{2}}\|\nabla\tilde{b}\|_2^{\frac{1}{2}}+\|\partial_z\tilde{b}
\|_2^{\frac{1}{2}}\|\nabla\partial_z\tilde{b}\|_2^{\frac{1}{2}}\Big)\|\nabla_H \tilde{b}\|_2^{\frac{1}{2}}\|\Delta\tilde{b}\|_2^{\frac{1}{2}}(\|\Delta \U\|_2+\|\partial_t \U\|_2)\notag\\
\leq&\frac{1}{8}\|\partial_t \U\|_2^2+\frac{1}{12}\|\Delta \U\|_2^2+\frac{1}{12}\|\Delta \B\|_2^2+C (\|\partial_z \tilde{b}\|_2^2\|\nabla\partial_z \tilde{b}\|_2^2+\|\tilde{b}\|_2^2\|\nabla\tilde{b}\|_2^2)\|\nabla\tilde{b}\|_2^2.
\end{align}

Thanks to the fact that $b_3$ is odd in $z$ and  $b_3|_{z=0}$, one gets

\begin{align}\label{c1}
b_3(x,y,z,t)=\int_0^z \partial_z b_3(x,y,\xi,t)d\xi=-\int_0^z \nabla_H\cdot \tilde{b}(x,y,\xi,t)d\xi,
\end{align}
and then, using  (\ref{c1}), we have
 \begin{align}\label{a6}
M_{32}=&\int_\Omega b_3\partial_z\tilde{b}(\partial_t \U-\Delta \U)\,dx dy dz\notag\\
\leq &\int_M\Big(\int_{-1}^1|\nabla_H\cdot\tilde{b}|dz\Big)\Big(\int_{-1}^1|\partial_z \tilde{b}|(|\Delta \U|+|\partial_t \U|)\,dz\Big)\,dxdy\notag\\
\leq& C\|\nabla_H\tilde{b}\|_2^{\frac{1}{2}}(\|\nabla_H\tilde{b}\|_2^{\frac{1}{2}}+\|\nabla_H^2\tilde{b}\|_2^{\frac{1}{2}})\|\partial_z\tilde{b}\|_2^\frac{1}{2}\notag\\
&\times(\|\partial_z \tilde{b}\|_2+\|\na_H\partial_z \tilde{b}\|_2)^{\frac{1}{2}}(\|\Delta \U\|_2+\|\partial_t \U\|_2)\notag\\
\leq&\|\nabla_H\tilde{b}\|_2^{\frac{1}{2}}\|\Delta_H\tilde{b}\|_2^{\frac{1}{2}}\|\partial_z\tilde{b}\|_2^{\frac{1}{2}}\|\nabla\partial_z\tilde{b}\|_2^{\frac{1}{2}}(\|\Delta \U\|_2+\|\partial_t \U\|_2)\notag\\
\leq&\frac{1}{8}\|\partial_t \U\|_2^2+\frac{1}{12}\|\Delta \U\|_2^2+\frac{1}{12}\|\Delta \tilde{b}\|_2+C\|\partial_z\tilde{b}\|_2^2\|\nabla\partial_z\tilde{b}\|_2^2\|\nabla\tilde{b}\|_2^2.
\end{align}

Next we turn to the next  term $M_4$
\begin{align*}
M_4=&\int_\Omega[(\tilde{b}\cdot\nabla_H) \U+b_3\partial_z \U](\Delta\tilde{b}-\partial_t\tilde{b})dxdydz\\
=&M_{41}+M_{42},
\end{align*}
again using the Lemma \ref{b7}, Poincar\'e  and Young inequalities
\begin{align}\label{a7}
M_{41}=&\int_{\Omega}( \tilde{b}\cdot\nabla_H)\U\cdot(\Delta \U-\partial_t \U)\,dxdydz\notag\\
\leq&\int_M\Big(\int_1^1(|\tilde{b}|+|\partial_z\tilde{b}|)dz\Big)\Big(\int_{-1}^1|\nabla_H \U|(|\Delta\tilde{b} |+|\partial_t \tilde{b}|)\,dz\Big)\,dxdy\notag\\
\leq&C\Big[\|\tilde{b}\|_2^{\frac{1}{2}}(\|\tilde{b}\|_2+\|\nabla_H\tilde{b}\|_2)^{\frac{1}{2}}+
\|\partial_z\tilde{b}\|_2^{\frac{1}{2}}(\|\partial_z\tilde{b}\|_2+\|\nabla_H\partial_z\tilde{b}\|_2)
^{\frac{1}{2}}\Big]\notag\\
&\times\|\nabla_H \U\|_2^{\frac{1}{2}}(\|\nabla_H \U\|_2+\|\nabla_H^2 \U\|_2)^{\frac{1}{2}}(\|\Delta \tilde{b}\|_2+\|\partial_t \tilde{b}\|_2)\notag\\
\leq&C\Big(\|\tilde{b}\|_2^{\frac{1}{2}}\|\nabla\tilde{b}\|_2^{\frac{1}{2}}+\|\partial_z\tilde{b}
\|_2^{\frac{1}{2}}\|\nabla\partial_z\tilde{b}\|_2^{\frac{1}{2}}\Big)\|\nabla_H \U\|_2^{\frac{1}{2}}\|\Delta_H \U\|_2^{\frac{1}{2}}(\|\Delta\tilde{b}\|_2+\|\partial_t \tilde{b} \|_2)\notag\\
\leq&\frac{1}{8}\|\partial_t \tilde{b}\|_2^2+\frac{1}{12}\|\Delta \U\|_2^2+\frac{1}{12}\|\Delta \tilde{b}\|_2^2+C (\|\partial_z \tilde{b}\|_2^2\|\nabla\partial_z \tilde{b}\|_2^2+\|\tilde{b}\|_2^2\|\nabla\tilde{b}\|_2^2)\|\nabla \U\|_2^2,
\end{align}
the second term  $M_{42}$ can also be bounded via Lemma \ref{b7}, Poincar\'e and Young inequalities
\begin{align}\label{a8}
M_{42}=&\int_\Omega b_3\partial_z \U(\partial_t\tilde{b} -\Delta \tilde{b})\,dx dy dz\notag\\
\leq &\int_M\Big(\int_{-1}^1|\nabla_H\cdot\tilde{b}|dz\Big)\Big(\int_{-1}^1|\partial_z \U|(|\Delta \tilde{b}|+|\partial_t \tilde{b}|)dz\Big)\,dxdy\notag\\
\leq& C\|\nabla_H\tilde{b}\|_2^{\frac{1}{2}}(\|\nabla_H\tilde{b}\|_2^{\frac{1}{2}}+\|\nabla_H^2
\tilde{b}\|_2^{\frac{1}{2}})\|\partial_z \tilde{u}\|_2^\frac{1}{2}\notag\\
&\times(\|\partial_z \U\|_2+\|\nabla_H\partial_z \U\|_2)^{\frac{1}{2}}(\|\Delta \tilde{b} \|_2+\|\partial_t \tilde{b}\|_2)\notag\\
\leq&\|\nabla_H\tilde{b}\|_2^{\frac{1}{2}}\|\Delta_H\tilde{b}\|_2^{\frac{1}{2}}\|\partial_z \U\|_2^{\frac{1}{2}}\|\nabla\partial_z \U\|_2^{\frac{1}{2}}(\|\Delta \tilde{b}\|_2+\|\partial_t \tilde{b}\|_2)\notag\\
\leq&\|\nabla\tilde{b}\|_2^{\frac{1}{2}}\|\Delta_H\tilde{b}\|_2^{\frac{1}{2}}\|\partial_z\U\|_2^{\frac{1}{2}}
\|\nabla\partial_z\U\|_2^{\frac{1}{2}}(\|\Delta \tilde{b}\|_2+\|\partial_t \tilde{b}\|_2)\notag\\
\leq&\frac{1}{8}\|\partial_t \tilde{b}\|_2^2+\frac{1}{12}\|\Delta \tilde{b}\|_2^2+C\|\nabla\tilde{b}\|_2^2\|\partial_z\U\|_2^2\|\nabla\partial_z\U\|_2^2.
\end{align}

Substituting these inequalities (\ref{a1})-(\ref{a8}) into (\ref{50}), we finally obtain
\begin{align*}
&\frac{d}{dt}(\|\nabla \U\|_2^2+\|\nabla\tilde{b}\|_2^2)+\frac{1}{2}(\|\Delta \U\|_2^2+\|\partial_t\U\|_2^2+\|\Delta\tilde{b} \|_2^2+\|\partial_t \tilde{b}\|_2^2)\\
\leq &C(\|\U\|_2^2\|\nabla \U\|_2^2+\|\partial_z \U\|_2^2\|\nabla\partial_z \U\|_2^2+\|\partial_z \tilde{b}\|_2^2\|\nabla\partial_z \tilde{b}\|_2^2+\|\tilde{b}\|_2^2\|\nabla\tilde{b}\|_2^2)(\|\nabla \U\|_2^2+\|\nabla \tilde{b}\|_2^2),
\end{align*}
recalling Proposition  \ref{d1} and Proposition \ref{d2}, the standard Gronwall inequality leads to
\begin{align*}
&\sup\limits_{0\leq s\leq t} \|\na\U\|_2^2(s)+\sup\limits_{0\leq s\leq t}\|\na\tilde{b}\|_2^2(s)+\frac{1}{2}\int_0^t(\|\De\U\|_2^2+\|\p_t\U\|_2^2+\|\De\B\|_2^2+\|\p_t\B\|_2^2)\,ds\\
\leq&exp\Big\{C\int_0^t(\|\U\|_2^2\|\na\U\|_2^2+\|\p_z\U\|_2^2\|\na\p_z\U\|_2^2+
\|\p_z\B\|_2^2\|\na\p_z\B\|_2^2+\|\B\|_2^2\|\na\B\|_2^2)\,ds\Big\}\\
&\times(\|\na\U\|_2^2+\|\na\B\|_2^2)\\
\leq&H(0),
\end{align*}
where
$$H(0)=(\|\na\U_0\|_2^2+\|\na\B_0\|_2^2)exp\big\{C(\|\U_0\|_2^4+\|\B_0\|_2^4+K^2(0))\big\},$$
which implies the conclusion.
\end{proof}
\vskip.1in
In this subsection, we deal with the second-order energy estimates, which are described by the following proposition.
\begin{prop}\label{p20}
(Second-order energy estimates). Suppose that $(\U_0,\tilde{b_0})\in H^2(\Omega)$. Then, we have the following estimates
\begin{align}\label{p21}
&\sup\limits_{0\leq s\leq t} \|\Delta \U\|_2^2(s)+\sup\limits_{0\leq s\leq t}\|\Delta\tilde{b}\|_2^2(s)+\frac{1}{2}\int_0^t(\|\na\Delta \U\|_2^2+\|\na\partial_t \U\|_2^2)ds\notag\\&+\frac{1}{2}\int_0^t(\|\na\Delta\tilde{b}\|_2^2+\|\na\partial_t\tilde{b}\|_2^2)\,ds
\leq exp\{CH^2(0)\}(\|\De\U_0\|_2^2+\|\De\B_0\|_2^2).
\end{align}
\end{prop}
\begin{proof}[Proof ] Taking the $L^2(\Om)$ inner product to the first equation and the third  equation of  (\ref{b5}) with $\De(\De\U-\p_t\U)$ and $\De(\De\B-\p_t\B)$, respectively,  and integration by parts,
\begin{align}\label{p22}
&\ddt\|\De\U\|_2^2+\ddt\|\De\B\|_2^2+\|\na\De\U\|_2^2+\|\na\p_t\U\|_2^2+\|\na\De\B\|_2^2+
\|\na\p_t\B\|_2^2\notag\\
=&\int_\Om\na[(u\C\na)\U]:\na(\De\U-\p_t\U)\,dxdydz-\int_\Om\na[(b\C\na)\B]:\na(\De\U-\p_t\U)\,dxdydz\notag\\
&+\int_\Om\na[(u\C\na)\B]:\na(\De\B-\p_t\B)\,dxdydz-\int_\Om\na[(b\C\na)\U]:\na(\De\B-\p_t\B)\,dxdydz\notag\\
=&N_1+N_2+N_3+N_4,
\end{align}
where $(:)$ denotes multiplication of two matrices.

To estimate  $N_1$, we use the Lemma \ref{b8}, the Poincar\'e and Young inequalities to get
\begin{align}\label{p23}
N_1\leq&\Big|\int_\Om\na[(u\C\na)\U]:\na(\De\U-\p_t\U)\,dxdydz\Big|\notag\\
\leq&\int_\Om[(\p_i u\C\na)\U+(u\C\p_i\na)\U]\C\p_i(\De\U-\p_t\U)\,dxdydz\notag\\
\leq&C(\|\p_i\na\U\|_2^{\frac{1}{2}}\|\p_i\De\U\|_2^{\frac{1}{2}}
\|\na\U\|_2^{\frac{1}{2}}\|\De\U\|_2^{\frac{1}{2}}+\|\na\U\|_2^{\frac{1}{2}}\|
\De\U\|_2^{\frac{1}{2}}\notag\\
&\|\p_i\na\U\|_2^{\frac{1}{2}}\|\p_i\De\U\|_2^{\frac{1}{2}})(\|\p_i\De\U\|_2+\|\p_i\p_t\U\|_2)\notag\\
\leq&\frac{1}{8}\|\na\De\U\|_2^2+\frac{1}{4}\|\na\p_t\U\|_2^2+C\|\na\U\|_2^2\|\De\U\|_2^4,
\end{align}
similar to the argument in (\ref{p23}), we have
\begin{align}\label{p24}
N_2\leq&\Big|\int_\Om\na[(b\C\na)\B]:\na(\De\U-\p_t\U)\,dxdydz\Big|\notag\\
\leq&\int_\Om[(\p_i b\C\na)\B+(b\C\p_i\na)\B]\C\p_i(\De\U-\p_t\U)\,dxdydz\notag\\
\leq&C(\|\p_i\na\B\|_2^{\frac{1}{2}}\|\p_i\De\B\|_2^{\frac{1}{2}}
\|\na\B\|_2^{\frac{1}{2}}\|\De\B\|_2^{\frac{1}{2}}+\|\na\B\|_2^{\frac{1}{2}}
\|\De\B\|_2^{\frac{1}{2}}\notag\\
&\|\p_i\na\B\|_2^{\frac{1}{2}}\|\p_i\De\B\|_2^{\frac{1}{2}})(\|\p_i\De\U\|_2+\|\p_i\p_t\U\|_2)\notag\\
\leq&\frac{1}{6}\|\na\De\B\|_2^2+\frac{1}{8}\|\na\De\U\|_2^2+\frac{1}{4}\|\na\p_t\U\|_2^2+C\|\na\B\|
_2^2\|\De\B\|_2^4,
\end{align}
and
\begin{align}\label{p24}
N_3\leq&\Big|\int_\Om\na[(u\C\na)\B]:\na(\De\B-\p_t\B)\,dxdydz\Big|\notag\\
\leq&\int_\Om[(\p_i u\C\na)\B+(u\C\p_i\na)\B]\C\p_i(\De\B-\p_t\B)\,dxdydz\notag\\
\leq&C(\|\p_i\na\U\|_2^{\frac{1}{2}}\|\p_i\De\U\|_2^{\frac{1}{2}}\|\na\B\|_2^{\frac{1}{2}}
\|\De\B\|_2^{\frac{1}{2}}+\|\na\U\|_2^{\frac{1}{2}}\|\De\U\|_2^{\frac{1}{2}}\notag\\
&\|\p_i\na\B\|_2^{\frac{1}{2}}\|\p_i\De\B\|_2^{\frac{1}{2}})(\|\p_i\De\B\|_2+\|\p_i\p_t\B\|_2)\notag\\
\leq&C(\|\De\U\|_2^{\frac{1}{2}}\|\na\De\U\|_2^\frac{1}{2}\|\na\B\|_2^{\frac{1}{2}}
\|\De\B\|_2^{\frac{1}{2}}+\|\na\U\|_2^{\frac{1}{2}}\|\De\U\|_2^{\frac{1}{2}}\notag\\
&\|\De\B\|_2^{\frac{1}{2}}\|\na\De\B\|_2^{\frac{1}{2}})(\|\na\De\B\|_2+\|\na\p_t\B\|_2)\notag\\
\leq&\frac{1}{4}\|\na\p_t\B\|_2^2+\frac{1}{8}\|\na\De\U\|_2^2+\frac{1}{6}\|\na\De\B\|_2^2\notag\\
&+C\|\De\U\|_2^2\|\na\B\|_2^2\|\De\B\|_2^2+C\|\na\U\|_2^2\|\De\U\|_2^2\|\De\B\|_2^2.
\end{align}

The last term $N_4$, we get
\begin{align}\label{p25}
N_4\leq&\Big|\int_\Om\na[(b\C\na)\U]:\na(\De\B-\p_t\B)\,dxdydz\Big|\notag\\
\leq&\int_\Om[(\p_i b\C\na)\U+(b\C\p_i\na)\U]\C\p_i(\De\B-\p_t\B)\,dxdydz\notag\\
\leq&C(\|\p_i\na\B\|_2^{\frac{1}{2}}\|\p_i\De\B\|_2^{\frac{1}{2}}
\|\na\U\|_2^{\frac{1}{2}}\|\De\U\|_2^{\frac{1}{2}}+\|\na\B\|_2^{\frac{1}{2}}
\|\De\B\|_2^{\frac{1}{2}}\notag\\
&\|\p_i\na\U\|_2^{\frac{1}{2}}\|\p_i\De\U\|_2^{\frac{1}{2}})(\|\p_i\De\B\|_2+\|\p_i\p_t\B\|_2)\notag\\
\leq&\frac{1}{6}\|\na\De\B\|_2^2+\frac{1}{8}\|\na\De\U\|_2^2+\frac{1}{4}\|\na\p_t\B\|_2^2\notag\\
&+C\|\De\B\|_2^2\|\De\U\|_2^2\|\na\U\|_2^2+C\|\na\B\|_2^2\|\De\B\|_2^2\|\De\U\|_2^2.
\end{align}

Collecting (\ref{p23})-(\ref{p25}) into (\ref{p22}), yields
\begin{align}\label{p26}
&\ddt(\|\De\U\|_2+\|\De\B\|_2^2)+\frac{1}{2}(\|\na\De\U\|_2^2+\|\na\p_t\U\|_2+\|\na\De\B\|_2^2+
\|\na\p_t\B\|_2^2)\notag\\
\leq&C(\|\na\U\|_2^2\|\De\U\|_2^2+\|\na\B\|_2^2\|\De\B\|_2^2)(\|\De\U\|_2^2+\|\De\B\|_2^2).
\end{align}

Applying the Gronwall inequality to the above inequalities, and recalling the Proposition \ref{b9}, one has
\begin{align*}
&\sup\limits_{0\leq s\leq t} \|\De\U\|_2^2(s)+\sup\limits_{0\leq s\leq t}\|\De\tilde{b}\|_2^2(s)+\frac{1}{2}\int_0^t(\|\na\De\U\|_2^2+\|\na\p_t\U\|_2^2+\|\na\De\B\|_2^2
+\|\na\p_t\B\|_2^2)\,ds\\
\leq&exp\Big\{C\int_0^t(\|\na\U\|_2^2\|\De\U\|_2^2+\|\na\B\|_2^2\|\De\B\|_2^2)\,ds\Big\}
(\|\De\U_0\|_2^2+\|\De\B_0\|_2^2)\\
\leq&exp\{CH^2(0)\}(\|\Delta\U_0\|_2^2+\|\De\B_0\|_2^2),
\end{align*}
where
$$H(0)=(\|\na\U_0\|_2^2+\|\na\B_0\|_2^2)exp\big\{C(\|\U_0\|_2^4+\|\B_0\|_2^4+K^2(0))\big\},$$
which implies the conclusion.
\end{proof}
 The global existence of strong solutions to the systems (\ref{b5}) is direct corollary of Proposition \ref{p20}.

\section{Uniqueness and continuous dependence of the strong solution}
\vskip .1in
In this section we will show the continuous dependence on the initial data and the uniqueness of the strong solution.
\vskip .1in
 Let $(\U_1, \B_1)$ and $(\U_2, \B_2)$  be two strong solutions of the system (\ref{b5}) with corresponding pressures $p_1$ and $p_2$, and the initial data $(\U_{0,1}, \B_{0,1})$ and $(\U_{0,2}, \B_{0,2})$,  respectively. Then the difference $X=\U_1-\U_2$, $q=p_1-p_2$ and $Y=\B_1-\B_2$ satisfy
\begin{align}\label{t5}
&\p_tX-\De X+\na_H q+(\U_1\cdot\na_H)X+(X\cdot\na_H)\U_2-\Big(\int_0^z \na_H\cdot\U_1(x,y,\xi,t)\,d\xi\Big)\p_zX\notag\\
&-\Big(\int_0^z\na_H\cdot X(x,y,\xi,t)\,d\xi\Big)\p_z\U_2-(\B_1\cdot\na_H)Y-(Y\cdot\na_H)\B_2\notag\\
&+\Big(\int_0^z\na_H\cdot\B_1(x,y,\xi,t)\,d\xi\Big)\p_z Y+\Big(\int_0^z\na_H\cdot Y(x,y,\xi,t)\,d\xi\Big)\p_z\B_2=0,
\end{align}
and
\begin{align}\label{t6}
&\p_tY+(\U_1\cdot\na_H)Y+(X\cdot\na_H)\B_2-\Big(\int_0^z \na_H\cdot\U_1(x,y,\xi,t)\,d\xi\Big)\p_zY\notag\\
&-\Big(\int_0^z\na_H\cdot X(x,y,\xi,t)\,d\xi\Big)\p_z\B_2-\De Y-(\B_1\cdot\na_H)X-(Y\cdot\na_H)\U_2\notag\\
&+\Big(\int_0^z\na_H\cdot\B_1(x,y,\xi,t)\,d\xi\Big)\p_z X
+\Big(\int_0^z\na_H\cdot Y(x,y,\xi,t)\,d\xi\Big)\p_z\U_2=0.
\end{align}

We equip the equations (\ref{t5}) and (\ref{t6}) with the following initial conditions,
\begin{align}\label{t7}
X(x,y,z,0)=\U_{0,1}-\U_{0,2},
\end{align}
and
\begin{align}\label{t8}
Y(x,y,z,0)=\B_{0,1}-\B_{0,2}.
\end{align}

\vskip .1in
 Taking the inner product of the  equation (\ref{t5}) with $X$ in $L^2(\Om)$, and the equation (\ref{t6}) with $Y$ in $L^2(\Om)$, we get
\begin{align*}
\frac{1}{2}\frac{d}{dt}\|X\|_2^2+\|\na X\|_2^2=&-\int_\Om\Big[(\U_1\cdot\na_H)X+(X\cdot\na_H)\U_2-\Big(\int_0^z \na_H\cdot\U_1(x,y,\xi,t)\,d\xi\Big)\p_zX\\
&-\Big(\int_0^z\na_H\cdot X(x,y,\xi,t)\,d\xi\Big)\p_z\U_2-(\B_1\cdot\na_H)Y-(Y\cdot\na_H)\B_2\\
&+\na_H q+\Big(\int_0^z\na_H\cdot\B_1(x,y,\xi,t)\,d\xi\Big)\p_z Y\\
&+\Big(\int_0^z\na_H\cdot Y(x,y,\xi,t)\,d\xi\Big)\p_z\B_2\Big]\cdot X\,dxdydz,
\end{align*}
and
\begin{align*}
\frac{1}{2}\frac{d}{dt}\|Y\|_2^2+\|\na Y\|_2^2=&-\int_\Om\Big[(\U_1\cdot\na_H)Y+(X\cdot\na_H)\B_2-\Big(\int_0^z \na_H\cdot\U_1(x,y,\xi,t)\,d\xi\Big)\p_zY\\
&-\Big(\int_0^z\na_H\cdot X(x,y,\xi,t)\,d\xi\Big)\p_z\B_2-(\B_1\cdot\na_H)X-(Y\cdot\na_H)\U_2\\
&+\Big(\int_0^z\na_H\cdot\B_1(x,y,\xi,t)\,d\xi\Big)\p_z X\\
&+\Big(\int_0^z\na_H\cdot Y(x,y,\xi,t)\,d\xi\Big)\p_z\U_2\Big]\cdot Y\,dxdydz.
\end{align*}

Applying integration by parts, and the boundary conditions (\ref{b1})-(\ref{b3}), we get
\begin{align}\label{t9}
-\int_\Om\Big((\U_1\cdot\na_H)X-\Big(\int_0^z\na_H\cdot\U_1(x,y,\xi,t)\,d\xi\Big)\p_z X\Big)\cdot X\,dxdydz=0,
\end{align}
\begin{align}\label{t10}
-\int_\Om\Big((\U_1\cdot\na_H)Y-\Big(\int_0^z\na_H\cdot\U_1(x,y,\xi,t)\,d\xi\Big)\p_z Y\Big)\cdot Y\,dxdydz=0,
\end{align}
and
\begin{align}\label{t11}
&\int_\Om\Big((\B_1\cdot\na_H)Y-\Big(\int_0^z\na_H\cdot\B_1(x,y,\xi,t)\,d\xi\Big)\p_z Y\Big)\cdot X\,dxdydz\notag\\
&+\int_\Om\Big((\B_1\cdot\na_H)X-\Big(\int_0^z\na_H\cdot\B_1(x,y,\xi,t)\,d\xi\Big)\p_z X\Big)\cdot Y\,dxdydz=0.
\end{align}

 Combining  (\ref{t9})- (\ref{t11}), we have
\begin{align*}
&\frac{1}{2}\frac{d}{dt}\|X\|_2^2+\|\na X\|_2^2+\frac{1}{2}\frac{d}{dt}\|Y\|_2^2+\|\na Y\|_2^2\\
=&-\int_\Om\Big[(X\cdot\na_H)\U_2
-(Y\cdot\na_H)\B_2-\Big(\int_0^z\na_H\cdot X(x,y,\xi,t)\,d\xi\Big)\p_z\U_2\notag\\
&+\Big(\int_0^z\na_H\cdot Y(x,y,\xi,t)\,d\xi\Big)\p_z \B_2\Big]\cdot X\,dxdydz
-\int_\Om\Big[(X\cdot\na_H)\B_2-(Y\cdot\na_H)\U_2\notag\\
&-\Big(\int_0^z\na_H\cdot X(x,y,\xi,t)\,d\xi\Big)\p_z\B_2
+\Big(\int_0^z\na_H\cdot Y(x,y,\xi,t)\,d\xi\Big)\p_z \U_2\Big]\cdot Y\,dxdydz.
\end{align*}

It follows from Lemma \ref{4} that
\begin{align}\label{t12}
\Big|\int_\Om(X\cdot\na_H)\U_2\cdot X\,dxdydz\Big|\leq&\|\na_H\U_2\|_2\|X\|_3\|X\|_6\notag\\
\leq &C\|\na\U_2\|_2\|X\|_2^{\frac{1}{2}}\|\na X\|_2^{\frac{3}{2}},
\end{align}
\begin{align}\label{t13}
\Big|\int_\Om(Y\cdot\na_H)\B_2\cdot X\,dxdydz\Big|\leq&\|\na_H\B_2\|_2\|Y\|_3\|X\|_6\notag\\
\leq & C\|\na\B_2\|_2\|\na X\|_2\|Y\|_2^{\frac{1}{2}}\|\na Y\|_2^{\frac{1}{2}},
\end{align}
\begin{align}\label{t14}
\Big|\int_\Om(X\cdot\na_H)\B_2\cdot Y\,dxdydz\Big|\leq&\|\na_H\B_2\|_2\|Y\|_3\|X\|_6\notag\\
\leq& C\|\na\B_2\|_2\|\na X\|_2\|Y\|_2^{\frac{1}{2}}\|\na Y\|_2^{\frac{1}{2}},
\end{align}
and
\begin{align}\label{t15}
\Big|\int_\Om(Y\cdot\na_H)\U_2\cdot Y\,dxdydz\Big|\leq&\|\na_H\U_2\|_2\|Y\|_3\|Y\|_6\notag\\
\leq& C\|\na\U_2\|_2\|Y\|_2^{\frac{1}{2}}\|\na Y\|_2^{\frac{3}{2}}.
\end{align}

Moreover
\begin{align}\label{19}
&\Big|\int_\Om\int_0^z\na_H\cdot X(x,y,\xi,t)\,d\xi\p_z\U_2\cdot X\,dxdydz\Big|\notag\\
\leq&\int_M\Big(\int_{-1}^{1}|\na_HX|\,dz\int_{-1}^{1}|\p_z\U_2||X|\,dz\Big)\,dxdy\notag\\
\leq&\int_M\Big(\int_{-1}^{1}|\na_H X|\,dz\Big(\int_{-1}^{1}|\p_z\U_2|^2\,dz\Big)^{\frac{1}{2}}\Big(\int_{-1}^{1}|X|^2\,dz\Big)^{\frac{1}{2}}\Big)\,dxdy\notag\\
\leq&\Big(\int_{M}\Big(\int_{-1}^{1}|\na X|\,dz\Big)^2\,dxdy\Big)^{\frac{1}{2}}\notag\\
&\times\Big(\int_M\Big(\int_{-1}^{1}|\p_z\U_2|^2\,dz\Big)^2\,dxdy\Big)^{\frac{1}{4}}
\Big(\int_M\Big(\int_{-1}^{1}|X|^2\,dz\Big)^2\,dxdy\Big)^{\frac{1}{4}}.
\end{align}

Using Cauchy-Schwarz inequality, one gets
\begin{align}\label{t16}
\Big(\int_M\Big(\int_{-1}^{1}|\na X|\,dz\Big)^2\,dxdy\Big)^{\frac{1}{2}}\leq C\|\na X\|_2.
\end{align}

Applying Sobolev inequality  (\ref{t24}) and Minkowsky inequality (\ref{t27}), one have
\begin{align}\label{t17}
\Big(\int_M\Big(\int_{-1}^{1}|X|\,dz\Big)^2\,dxdy\Big)^{\frac{1}{2}}\leq &C\int_{-1}^{1}\Big(\int_M|X|^4\,dxdy\Big)^{\frac{1}{2}}\,dz\notag\\
\leq& C\int_{-1}^1|X||\na X|\,dz\notag\\
\leq& C\|X\|_2\|\na X\|_2,
\end{align}
and
\begin{align}\label{t18}
\Big(\int_M\Big(\int_{-1}^{1}|\p_z\U_2|\,dz\Big)^2\,dxdy\Big)^{\frac{1}{2}}\leq &C\int_{-1}^{1}\Big(\int_M|\p_z\U_2|^4\,dxdy\Big)^{\frac{1}{2}}\,dz\notag\\
\leq &C\int_{-1}^1|\p_z\U_2||\na \p_z\U_2|\,dz\notag\\
\leq& C\|\p_z\U_2\|_2\|\na \p_z\U_2\|_2.
\end{align}

Therefore, by estimates (\ref{t16})-(\ref{t18}), we get
\begin{align}\label{t20}
&\Big|\int_\Om\int_0^z\na_H\cdot X(x,y,\xi,t)\,d\xi\p_z\U_2\cdot X\,dxdyz\Big|\notag\\
\leq&C\|\na X\|_2^{\frac{3}{2}}\|X\|_2^{\frac{1}{2}}\|\p_z\U_2\|_2^{\frac{1}{2}}\|\na\p_z\U_2\|_2^{\frac{1}{2}}.
\end{align}

It follows from the similar argument in (\ref{t20}) that
\begin{align}\label{t21}
&\Big|\int_\Om\int_{0}^{z}\na_H\cdot Y(x,y,\xi,t)\,d\xi\p_z\B_2\C X\,dxdydz\Big|\notag\\
\leq&\|\na Y\|_2\|X\|_2^{\frac{1}{2}}\|\na X\|_2^{\frac{1}{2}}\|\p_z\B_2\|_2^{\frac{1}{2}}\|\na\p_z\B_2\|_2^{\frac{1}{2}},
\end{align}
and
\begin{align}\label{t22}
&\Big|\int_\Om\int_{0}^{z}\na_H\cdot X(x,y,\xi,t)\,d\xi\p_z\B_2\C Y\,dxdydz\Big|\notag\\
\leq&\|\na X\|_2\|\p_z\B_2\|_2^{\frac{1}{2}}\|\na \p_z\B_2\|_2^{\frac{1}{2}}\|Y\|_2^{\frac{1}{2}}\|\na Y\|_2^{\frac{1}{2}},
\end{align}
moreover,
\begin{align}\label{t23}
&\Big|\int_\Om\int_{0}^{z}\na_H\cdot Y(x,y,\xi,t)\,d\xi\p_z\U_2\C Y\,dxdydz\Big|\notag\\
\leq&\|\p_z\U_2\|_2^{\frac{1}{2}}\|\na\p_z\U_2\|_2^{\frac{1}{2}}\|Y\|_2^{\frac{1}{2}}\|\na Y\|_2^{\frac{3}{2}}.
\end{align}

Therefore, substituting estimates (\ref{t12})-(\ref{t15}) and  (\ref{t20})-(\ref{t23}), we reach
\begin{align*}
&\frac{1}{2}\frac{d}{dt}(\|X\|_2^2+\|Y\|_2^2)+\|\na X\|_2^2+\|\na Y\|_2^2\\
\leq&C\big(\|\na\U_2\|_2+\|\p_z\U_2\|_2^{\frac{1}{2}}\|\na\p_z\U_2\|_2^{\frac{1}{2}}\big)\|X\|_2^{\frac{1}{2}}\|\na X\|_2^{\frac{3}{2}}\\
&+C\|\na\B_2\|_2\|\na X\|_2\|Y\|_2^{\frac{1}{2}}\|\na Y\|_2^{\frac{1}{2}}\\
&+C\|\na Y\|_2\|X\|_2^{\frac{1}{2}}\|\na X\|_2^{\frac{1}{2}}\|\p_z\B_2\|_2^{\frac{1}{2}}\|\na\p_z\B_2\|_2^{\frac{1}{2}}\\
&+C\|\na X\|_2\|Y\|_2^{\frac{1}{2}}\|\na Y\|_2^{\frac{1}{2}}\|\p_z\B_2\|_2^{\frac{1}{2}}\|\na\p_z\B_2\|_2^{\frac{1}{2}}\\
&+C\Big(\|\na\U_2\|_2+\|\p_z\U_2\|_2^{\frac{1}{2}}\|\na\p_z\U_2\|_2^{\frac{1}{2}}\Big)\|Y\|_2^{\frac{1}{2}}\|\na Y\|_2^{\frac{3}{2}}.
\end{align*}

It follows from Young inequality that
\begin{align*}
&\frac{1}{2}\frac{d}{dt}(\|X\|_2^2+\|Y\|_2^2)+\|\na X\|_2^2+\|\na Y\|_2^2\\
\leq&C\Big(\|\na\U_2\|_2^4+\|\na\B_2\|_2^4+\|\na\U_2\|_2^2+\|\p_z\U_2\|_2^2\|\na\p_z\U_2\|_2^2+\|\p_z\B_2\|_2^2\|\na\p_z\B_2\|_2^2
\Big)\\
&\times(\|X\|_2^2+\|Y\|_2^2).
\end{align*}

Thanks to Gronwall inequality,
\begin{align*}
\|X\|_2^2+\|Y\|_2^2
\leq &\big(\|X(t=0)\|_2^2+\|Y(t=0)\|_2^2\big)\\
&exp\Big\{C\int_0^t(\|\na\U_2(s)\|_2^4+\|\na\B_2(s)\|_2^4+\|\na\U_2(s)\|_2^2\\
&+\|\p_z\U_2(s)\|_2^2\|\na\p_z\U_2(s)\|_2^2+\|\p_z\B_2(s)\|_2^2\|\na\p_z\B_2(s)\|_2^2)\,ds\Big\}.
\end{align*}

Since $(\U_2,\B_2)$ is a strong solution,
\begin{align*}
\|X\|_2^2+\|Y\|_2^2
\leq &\big(\|X(t=0)\|_2^2+\|Y(t=0)\|_2^2\big)exp\{C(H^2(0)+H(0)+K(0)H(0))\}.
\end{align*}

The above inequality proves the continuous dependence of the solutions on the initial data. In particular, $X(t=0)=Y(t=0)=0$, we have $X(t)=Y(t)=0$, for all $t\geq 0$. Therefore, the strong solution is unique.

\vskip .2in
\centerline{\bf Acknowledgments}
\vskip .1in
This work is supported by the National Natural Science Foundation of China grant 11971331, 12125102, and Sichuan Youth  Science and Technology Foundation 2021 JDTD0024.
 \vskip .2in

{\bf Conflict of interest. } The authors declare that they have no
conflict of interest and our manuscript has no associated data.

\end{document}